\documentclass[11pt]{amsart}
\usepackage{mathrsfs}
\usepackage{amssymb}

\usepackage{amscd}
\usepackage{amsmath, amssymb}
\usepackage{amsfonts}
\usepackage[colorlinks,linkcolor=blue,citecolor=blue, pdfstartview=FitH]{hyperref}
\usepackage[all]{xy}

\numberwithin{equation}{section}

\theoremstyle{plain}
\newtheorem{thm}{Theorem}[section]
\newtheorem{theorem}[thm]{Theorem}
\newtheorem{lemma}[thm]{Lemma}
\newtheorem{corollary}[thm]{Corollary}
\newtheorem{proposition}[thm]{Proposition}

\theoremstyle{definition}

\newtheorem{remark}[thm]{Remark}

\newtheorem{definition}[thm]{Definition}

\newtheorem{example}[thm]{Example}

\newtheorem{defn-thm}[thm]{Definition-Theorem}

\newcommand{\Image}{{ \textrm{im}\,}}

\newcommand{\wg}{{\,\wedge\,}}
\begin{document}

\title{Deformations of Dolbeault cohomology classes}
\author{Wei Xia}
\address{Wei Xia, Mathematical Science Research Center, Chongqing University of Technology, Chongqing, P.R.China, 400054.} \email{xiawei@cqut.edu.cn, xiaweiwei3@126.com}

\thanks{This work was supported by the National Natural Science Foundation of China No. 11901590 and Scientific Research Foundation of Chongqing University of Technology.}
\date{\today}

\begin{abstract}
In this paper, we establish a deformation theory for Dolbeault cohomology classes valued in holomorphic tensor bundles. We prove the extension equation which will play the role of Maurer-Cartan equation. Following the classical theory of Kodaira-Spencer-Kuranishi, we construct a canonical complete family of deformations by using the power series method. We also prove a simple relation between the existence of deformations and the varying of the dimensions of Dolbeault cohomology. The deformations of $(p,q)$-forms is shown to be unobstructed under some mild conditions. By analyzing Nakamura's example of complex parallelizable manifolds, we will see that the deformation theory developed in this work provides precise explanations to the jumping phenomenon of Dolbeault cohomology.

\vskip10pt
\noindent
{\bf MSC~Classification (2020):} 32G05, 32L10, 55N30, 32G99
\end{abstract}
\maketitle

\tableofcontents

\section{Introduction}
The foundations of the theory of deformations of complex analytic structures on a compact complex manifold $X$ is developed around the 1960s' mainly by Kodaira and Spencer~\cite{KS58,KS60}. In \cite{KNS58}, Kodaira-Nirenberg-Spencer proved a fundamental existence theorem for the deformations of complex analytic structures under the assumption that $H^2(X, \Theta)=0$, where $\Theta$ is the sheaf of germs of holomorphic vector fields on $X$. Later on, Kuranishi~\cite{Kur62} proved the theorem of existence in full generality which states that there exists a semiuniversal family of deformations of $X$ for any compact complex manifold $X$. The family of complex manifolds $\pi: \mathcal{X}\to \mathcal{B}$ thus constructed is called the \emph{Kuranishi family}, and $\mathcal{B}$ is called the \emph{Kuranishi space} of $X$. Note that because we are studying small deformations, only the germs of $\pi$ and $\mathcal{B}$ are well-defined notions.

In \cite{Tod89}, Todorov constructed a canonical family of holomorphic $(n,0)$-forms for the Kuranishi family of Calabi-Yau manifolds.
Liu-Sun-Yau~\cite{LSY09} studied deformations of pluricanonical forms on the Kuranishi space in the case of Riemann surfaces. A remarkable feature of their study is the use of the extension operator and the iteration procedure to construct the deformations. By means of the extension operator, it is natural to consider the extension equation which plays the role of Maurer-Cartan equation. The iteration procedure is quite analogous to the power series method employed in the classical theory of deformations of complex structures~\cite{MK71,Kod86}. On the other hand, in Griffiths' visionary work~\cite{Gri64}, an abstract algebraic framework has been established for the extension problem in Dolbeault cohomology. The present work is much inspired by these ideas. See also \cite{LZ18,RZ18,RWZ19,LRY15,Sun12,Wav73} for some related works.

The main aim of this paper is to develop a deformation theory in the sense of Kodaira-Spencer for Dolbeault cohomology classes (valued in holomorphic tensor bundles) and consider some of its applications. As in the classical theory, the so-called Beltrami differential is a key player in our story. Let $X$ be a compact complex manifold, any small deformation $X_t$ of $X$ can be represented by a Beltrami differential $\phi = \phi(t) \in A^{0,1}(X, T^{1,0})$ (see e.g. \cite{Voi02I,Huy05,Kur65}), which may be viewed as a smooth bundle map $\phi : T^{0,1}\to T^{1,0}$, such that the Maurer-Cartan equation is satisfied:
\begin{equation}\label{Maurer-Cartan eq 0}
\bar{\partial}\phi - \frac{1}{2}[\phi ,\phi ]=0.
\end{equation}
By means of the Beltrami differential $\phi$, we can introduce the \emph{extension operator} $\rho: A^{0,q}(X,E)\to A^{0,q}(X_t,E_t)$ for $E$-valued $(0,q)$-forms, where $E$ and $E_t$ are corresponding holomorphic tensor bundles on $X$ and $X_t$, respectively. The extension operator $\rho$ adopted in this paper is a combination of the exponential operator studied in \cite{Tod89,LSY09,LRY15} and the projection operator used in \cite{Ham77,Hua95}. Given $\sigma\in A^{0,q}(X,E)$, the \emph{extension equation} tells us when is $\rho\sigma\in A^{0,q}(X_t,E_t)$ a $\bar{\partial}_t$-closed section where $\bar{\partial}_t$ is the Dolbeault operator on $X_t$. In fact, we have the following commutative diagram (see Section \ref{sec-Extension isomorphism}):
\[
\xymatrix{
  A^{0,\bullet}(X, E) \ar[d]_{\bar{\partial}_{\phi(t)} } \ar[r]^{\rho} & A^{0,\bullet}(X_t, E_t) \ar[d]^{\bar{\partial}_t} \\
  A^{0,\bullet+1}(X, E) \ar[r]^{\rho} & A^{0,\bullet+1}(X_t, E_t),   }
\]
which implies that
\begin{equation}\label{general extension equation 0}
\bar{\partial}_t\rho\sigma=0\Longleftrightarrow\bar{\partial}\sigma - \langle \phi | \sigma \rangle =0.
\end{equation}
As mentioned above, \eqref{general extension equation 0} will play the same role as the Maurer-Cartan equation~\eqref{Maurer-Cartan eq 0}. It should be pointed out that various special cases of \eqref{general extension equation 0} already appeared in previous works. In fact, if $E$ is trivial and $q=0$, i.e. for functions, \eqref{general extension equation 0} reduced to Proposition 1.2 of~\cite[pp.\,151]{MK71} which is a well-known fact. If $E$ is the tangent bundle $T_X^{1,0}$, this was proved by Hamilton~\cite[pp.\,7]{Ham77} and Zhang~\cite[pp.\,55]{Zha14} ($q=0$). For $q=0$, if $E$ is the canonical bundle $K_X$, this was proved by Liu-Rao-Yang~\cite[Prop.\,5.1]{LRY15}; if $E$ is the pluricanonical bundle $K_X^m$, this was proved by Liu-Zhu~\cite[Lem.\,6.1]{LZ18}; the case of $E=\Omega^p$ was also proved by Rao-Zhao~\cite[Prop.\,2.13]{RZ18}.

A useful consequence of the above commutative diagram is the following isomorphism induced by $\rho$,
\begin{equation}\label{ext-iso for coho 0}
H^{0,q}_{\bar{\partial}_{\phi}}(X,E) \cong H^{0,q}_{\bar{\partial}_{t}}(X_t,E_t),~\forall q \geq 0,
\end{equation}
where
\[
H^{0,q}_{\bar{\partial}_{\phi}}(X,E) := \frac{\ker \bar{\partial}_{\phi} \cap A^{0,q}(X,E)}{\Image \bar{\partial}_{\phi} \cap A^{0,q}(X,E)},~\forall q\geq 0,
\]
and
\[
\bar{\partial}_{\phi}:= \bar{\partial}-\langle\phi|~:~A^{0,\bullet}(X, E)\longrightarrow A^{0,\bullet+1}(X, E).
\]
The paring $\langle\cdot|\cdot\rangle : A^{0,k}(X, T_{X}^{1,0})\times A^{0,q}(X, E)  \rightarrow  A^{0,q+k}(X, E)$ and the operator $\langle\phi|$ are natural generalization of the Lie derivative and the Fr\"olicher-Nijenhuis bracket. It is commonly believed that every deformation problem is controlled by a DGLA~\cite{Ma05,SS12}, we find that in our case the deformation is determined by modules over the deformed Kodaira-Spencer DGLA $(A^{0,\bullet}(X, T_{X}^{1,0}), \bar{\partial}_{\phi}, [\cdot,\cdot])$. Motivated by these facts, we make the following definitions (see Section \ref{Universal property for the canonical deformation}).

Let $\pi: (\mathcal{X}, X)\to (B,0)$ be a small deformation~\footnote{Though the base space $B$ may be a singular complex space, such as the Kuranishi space, only deformations of complex manifolds will be considered. Thus every fiber of $\pi$ are complex manifolds. This is our basic setting in this paper.} of a compact complex manifold $X$ such that for each $t\in B$ the complex structure on $X_t$ is represented by Beltrami differential $\phi(t)$ and $E$ be a holomorphic tensor bundle on $X$. Given $y\in \ker\bar{\partial}\cap A^{0,q}(X,E)$ and $T\subseteq B$, which is a complex subspace of $B$ containing $0$, a \emph{deformation} of $y$ (w.r.t. $\pi: (\mathcal{X}, X)\to (B,0)$ ) on $T$ is a family of $E$-valued $(0,q)$-forms $\sigma (t)$ such that
\begin{itemize}
  \item[1.] $\sigma (t)$ is holomorphic in $t$ and $\sigma (0) = y$;
  \item[2.] $\bar{\partial}_{\phi(t)}\sigma (t) = \bar{\partial}\sigma (t) - \langle\phi (t) | \sigma (t) \rangle =0,~\forall t\in T$.
\end{itemize}
We say $y$ has \emph{unobstructed deformation w.r.t. $\pi$} if $B$ is smooth and there exists a deformation of $y$ on $B$. We say $y$ has \emph{unobstructed deformation} if for any small deformation $\pi: (\mathcal{X}, X)\to (B,0)$ of $X$ with smooth $B$ there exists a deformation of $y$ (w.r.t. $\pi$) on $B$.
Given a Dolbeault cohomology class $\alpha\in H_{\bar{\partial}}^{0,q}(X,E)$, a \emph{deformation} of $\alpha$ (w.r.t. $\pi$) on $T$ is defined as a triple $(y,\sigma (t),T)$ which consisting of a representative $y\in\alpha$ and a deformation $\sigma (t)$ of $y$ (w.r.t. $\pi$) on $T$. Two deformations (w.r.t. $\pi$) $(y,\sigma (t),T)$ and $(y',\sigma' (t),T)$ of $\alpha$ on $T$ are \emph{equivalent} if
\[
[\sigma (t) - \sigma' (t)] = 0 \in H^{0,q}_{\bar{\partial}_{\phi(t)}}(X,E),~\forall t\in T.
\]
We say $\alpha\in H_{\bar{\partial}}^{0,q}(X,E)$ has \emph{unobstructed deformation} if for any small deformation $\pi: (\mathcal{X}, X)\to (B,0)$ of $X$ with smooth $B$ there exists $y\in\alpha$ such that $y$ has unobstructed deformation w.r.t. $\pi$.
If for any small deformation $\pi: (\mathcal{X}, X)\to (B,0)$ of $X$ with smooth $B$, every Dolbeault cohomology classes in $H_{\bar{\partial}}^{0,q}(X,E)$ have unobstructed deformation w.r.t. $\pi$, then we say the \emph{deformations of classes in $H_{\bar{\partial}}^{0,q}(X,E)$ are unobstructed}.
Now, assume $X$ has been equipped with a fixed Hermitian metric, a deformation $\sigma (t)$ of $y\in \ker\bar{\partial}\cap A^{0,q}(X,E)$ on $T$ (w.r.t. $\pi$) is \emph{canonical} if it satisfies
$\sigma(t) = y + \bar{\partial}^*G \langle\phi(t) | \sigma(t)\rangle$ for any $t\in T$ where $G$ is the $\bar{\partial}$-Green's operator. A deformation $(y,\sigma (t),T)$ of $[y]\in H_{\bar{\partial}}^{0,q}(X,E)$ (w.r.t. $\pi$) is \emph{canonical} if $\sigma (t)$ is a canonical deformation of $y$.

As in the classical Kodaira-Spencer theory, a fundamental question is how to construct a complete/universal family of deformations. This is answered by the following result (see Theorem \ref{Deformation of Dolbeault cohomology classes: general case} and Theorem \ref{universal and uniqueness of the canonical deformation}):
\begin{theorem}\label{main result combined 0}
Let $\pi: (\mathcal{X}, X)\to (B,0)$ be a small deformation of a compact Hermitian manifold $X$ such that for each $t\in B$ the complex structure on $X_t$ is represented by a Beltrami differential $\phi(t)$. Suppose $E$ is a holomorphic tensor bundle on $X$ which is equipped with a fixed Hermitian metric and $\mathcal{H}^{0,q}(X,E)$ is the space of harmonic $E$-valued $(0,q)$-forms.
\begin{itemize}
	\item[$(i)$] For any linear subspace $V\subseteq\mathcal{H}^{0,q}(X,E)$, there is a maximal analytic subset $B(V)$ of $B$ such that the canonical deformation of $\sigma_0$ exists on $B(V)$ where $\sigma_0\in V$ is arbitrary. Furthermore, assume $S$ is an analytic subset of $B$ with $0\in S$ and $y\in \ker\bar{\partial}\cap A^{0,q}(X,E)$. Then the canonical deformation of $y$ exist on $B(\mathbb{C}\mathcal{H}y)$ and if the canonical deformation of $y$ exists on $S$ we must have $S\subseteq B(\mathbb{C}\mathcal{H}y)$;
	\item[$(ii)$] For any deformed Dolbeault cohomology class $[u]\in H_{\bar{\partial}_{\phi(t)}}^{0,q}(X,E)$, there exists $\sigma_0\in \mathcal{H}^{0,q}(X,E)$ such that $[u]=[\sigma(t)]$ where $\sigma(t)$ is the canonical deformation of $\sigma_0$;
    \item[$(iii)$] For any class $\alpha\in H_{\bar{\partial}}^{0,q}(X,E)$, if the canonical deformation of some representative $y\in\alpha$ exists on $S\subseteq B$ then the canonical deformation of any other representative $y'\in\alpha$ also exists on $S\subseteq B$ and their canonical deformations are equivalent.
    \item[$(iv)$]  Let $\varpi: (\mathcal{Y}, Y_{s_0}) \to (D, s_0)$ be a pullback of the Kuranishi family $\pi_K: (\mathcal{X},X)\to (\mathcal{B},0)$ with the following commutative diagram:
\[
\xymatrix{
(\mathcal{Y}, Y_{s_0})  \ar[r]^{\Phi} \ar[d]^{\varpi} & (\mathcal{X}, X) \ar[d]^{\pi_K} \\
   (D, s_0)  \ar[r]^{h}    & (\mathcal{B},0) , }
\]
and $\Phi_{s}:=\Phi\mid_{Y_{s}}$ for each $s\in D$. Assume $Y_{s_0}$ is equipped with the induced Hermitian metric so that $\Phi_{s_0}^*:Y_{s_0}\to X$ is an isometry and $\tau(s)$ is the canonical deformation of $\Phi_{s_0}^*y$ on $\tilde{D}\subseteq D$. Then $\tilde{D}\subseteq h^{-1}(\mathcal{B}(\mathbb{C}\mathcal{H}y))$ and $\tau(s)$ is just the pullback of the canonical deformation of $y$ on $\mathcal{B}(\mathbb{C}\mathcal{H}y)$.
\end{itemize}
\end{theorem}

The canonical deformation $\sigma(t)\in A^{0,q}(X,E)$ depends uniquely on the initial value $\sigma(0)=y$ and can be written out in an iteration form. It may happen that $\sigma(t)$ is $\bar{\partial}_{\phi(t)}$-exact for some $t\in B(\mathbb{C}\mathcal{H}y)$ if $q>0$. Along with the possible obstructions of the canonical deformations, these two facts\footnote{It turns out that the effect of the former is equivalent to the obstructions of classes in $H^{0,q-1}_{\bar{\partial}}(X,E)$, see Proposition \ref{prop-kerdbar*-imagedbarphi}.} are the exact reasons which cause $\dim H^{q}(X_t,E_t)$ to jump, see Section \ref{Examples of obstructed deformations and the jumping phenomena} for concrete examples.

A given Dolbeault cohomology class may have many inequivalent deformations (See the paragraph just before Remark~\ref{rk-deformation nonclosed}) and the relations between these inequivalent deformations is quite mysterious. This is why we insist on considering canonical deformations\footnote{Even for canonical deformations, it is not clear how they depend on the choices of the Hermitian metric. For example, is it true that the equivalence class of canonical deformations of a given class is unique (thus independent of the Hermitian metrics chosen)?} in this paper. Also, it is still not clear whether (and how) the existence of
deformations is related to the varying of $\dim H^{q}(X_t,E_t)$ in general. Nonetheless, we have the following:
\begin{theorem}[=Theorem \ref{thm-2nd-main}]\label{thm-2nd-main-0}
Let $\pi: (\mathcal{X}, X)\to (B,0)$ be a small deformation of a compact Hermitian manifold $X$ such that for each $t\in B$ the complex structure on $X_t$ is represented by a Beltrami differential $\phi(t)$. Suppose $E$ is a holomorphic tensor bundle on $X$ which is equipped with a fixed Hermitian metric. For each $q\geq 0$ and $t\in B$, set
\[
v^q_t:=\dim H_{\bar{\partial}}^{0,q}(X,E)-\dim \ker\bar{\partial}_{\phi(t)}\cap\ker\bar{\partial}^*\cap A^{0,q}(X,E) \geq 0,
\]
then we have $(v^{-1}_t:=0)$
\begin{equation}\label{eq-1461-0}
\dim H_{\bar{\partial}}^{0,q}(X,E)=\dim H_{\bar{\partial}_t}^{0,q}(X_t,E_t)+v^q_t+v^{q-1}_t.
\end{equation}
In particular, $\dim H_{\bar{\partial}_t}^{0,q}(X_t,E_t)$ is independent of $t\in B$ if and only if the canonical deformations of classes in $H_{\bar{\partial}}^{0,q}(X,E)$ and $H_{\bar{\partial}}^{0,q-1}(X,E)$ exist on $B$.
\end{theorem}

It is interesting to compare \eqref{eq-1461-0} with a classical result of Kodaira-Spencer (see formula $(43)$ in \cite[pp.\,68]{KS60}). In their formula,  $v^q_t$ is defined via the spectrum of Laplacian operators. Here our $v^q_t$ in \eqref{eq-1461-0} is determined by the properties of canonical deformations.
Moreover, Theorem \ref{thm-2nd-main-0} gives a solution to the question of Rao-Zhao \cite[Question\,1.6]{RZ18} where they ask for sufficient and necessary conditions for the deformation invariance of Hodge numbers.

A celebrated theorem proved by Siu~\cite{Siu02a} says that the plurigenera are deformation invariants for projective family of algebraic manifolds where projective family means that there is a positive line bundle on the total space of the family~\cite{Siu98}. It seems not known whether this is still true if we only assume the central fiber of the analytical family is projective. Nevertheless, it is conjectured by Siu~\cite{Siu02a,Siu02b} that the plurigenera is a deformation invariant for compact K\"ahler manifolds, see \cite{RT20,Dem20,CP20,Kol21} for some recent progress. Furthermore, Nakayama even suspect it is enough to assume the central fiber is in the Fujiki class $\mathcal{C}$~\cite{Nak04}. It is proposed by Liu-Zhu that we can approach Siu's conjecture by using the iteration method \cite[Sec.\,6]{LZ18}. Theorem \ref{thm-2nd-main-0} shows that in order to solve Siu's conjecture we only need to show the deformations of pluricanonical forms on a compact K\"ahler manifold is canonically unobstructed, see Section \ref{Universal property for the canonical deformation} for the definition of canonically unobstructedness.

For the deformations of $(p,q)$-forms, we have the following result:
\begin{theorem}[=Theorem \ref{thm-unobstructed-def for p,q-forms}]\label{thm-unobstructed-def for p,q-forms-0}
Let $X$ be a compact Hermitian manifold.
\begin{itemize}
  \item[$(i)$] Assume $\partial_{A,\bar{\partial}(\ker\partial)}^{p,q}=0$ and $\partial_{A,\bar{\partial}}^{p-1,q+1}=0$ where
\[
\partial_{A,\bar{\partial}(\ker\partial)}^{p,q}:H_{A}^{p,q}(X)\longrightarrow \frac{\ker\bar{\partial}\cap A^{p+1,q}(X)}{\bar{\partial}\left(\ker\partial\cap A^{p+1,q-1}(X) \right)},\quad \partial_{A,\bar{\partial}}^{p,q}:H_{A}^{p,q}(X)\longrightarrow H_{\bar{\partial}}^{p+1,q}(X),
\]
is the natural map induced by $\partial$ and $H_{A}^{p,q}(X)=\frac{\ker \partial\bar{\partial} \cap  A^{p,q}(X)}{ \partial A^{p-1,q}(X)+\bar{\partial}A^{p,q-1}(X) }$ is the Aeppli cohomology of $X$. Then the deformations of classes in $H_{\bar{\partial}}^{p,q}(X)$ are canonically unobstructed;
  \item[$(ii)$] Let $\pi: (\mathcal{X}, X)\to (B,0)$ be a small deformation of $X$ such that $B$ is smooth. Assume $\partial_{A,\bar{\partial}(\ker\partial)}^{p,r}=0$ and $\partial_{A,\bar{\partial}}^{p-1,r+1}=0$ for $r=q,q-1$, then $\dim H_{\bar{\partial}_t}^{p,q}(X_t)$ is independent of $t\in B$.
\end{itemize}
\end{theorem}
Recently, Rao-Zhao~\cite{RZ18} studied the deformation invariance of Hodge numbers by assuming a weak form of the $\partial\bar{\partial}$-lemma. It is proved in \cite[Thm.\,1.3]{RZ18} that if $H_{BC}^{p+1,q}(X)\to H_{\partial}^{p+1,q}(X)$ is injective, $\partial_{A,\bar{\partial}}^{p-1,q+1}=0$ and $\dim H_{\bar{\partial}_t}^{p,q-1}(X_t)$ is deformation invariant then $\dim H_{\bar{\partial}_t}^{p,q}(X_t)$ is also deformation invariant. This can be implied by Theorem \ref{thm-unobstructed-def for p,q-forms-0} because $\partial_{A,\bar{\partial}(\ker\partial)}^{p,q}=0$ is a consequence of the injectivity of $H_{BC}^{p+1,q}(X)\to H_{\partial}^{p+1,q}(X)$.

As to obstructed deformations, we will study the examples provided by Nakamura~\cite{Nak75} in Section~\ref{Examples of obstructed deformations and the jumping phenomena}. The theory developed in this paper gives a precise explanation to the jumping phenomenon, see~\cite{Ye08,Ye10,Wav73,Kle71} for related works on this topic. Note that the techniques employed in this article provides a new method to compute the Dolbeault cohomology for small deformations of compact complex manifolds. Compare~\cite{AK17b}.

\section{The extension operator and extension equation} \label{The extension operator and extension equation}
In this section, we introduce the definition of the extension operator and derive the extension equation for sections valued in a holomorphic tensor bundle which will be fundamental for our purpose.
\subsection{The extension operator}
\label{The extension operator}
Let $\pi : \mathcal{X} \to \Delta^k$ be an analytic family of compact complex manifolds (i.e. a proper surjective holomorphic submersion between complex manifolds) of dimension $n$ over the unit polydisc $\Delta^k$ in $\mathbb{C}^k$ with fiber $X_t = \pi^{-1}(t)$, where $t = (t_1,\cdots, t_k)\in \Delta^k$.
\subsubsection{Beltrami differentials}
By the Ehresmann theorem, we have the following commutative diagram
\[
  \xymatrix{
  X\times \Delta^k \ar[dr] \ar[r]^-{F}
                & \mathcal{X} \ar[d]^-{\pi}  \\
                & \Delta^k, }
\]
where $F$ is a diffeomorphism and $X:=X_0$. For any $t\in \Delta^k$ set $f_t:=F\mid_{X\times \{t\}}$, then $f_t$ is a diffeomorphism from $X$ to $X_t$. Each fiber $X_t$ is then regarded as a deformation (of complex structure) with respect to $X$. For sufficiently small $t$, a
deformation $X_t$ is represented by a Beltrami differential $\phi(t)\in A^{0,1}(X, T^{1,0})$.
Let $z^1,\cdots, z^n$ and $w^1,\cdots, w^n$ be holomorphic coordinates on $X$ and $X_t$ respectively, the Beltrami differential can be defined by
\begin{equation}\label{Beltrami differential coord}
\phi(t):=\phi^{\alpha}_{\beta} d\bar{z}^{\beta}\otimes\frac{\partial}{\partial z^{\alpha}}=\left(\frac{\partial w}{\partial z}\right)^{-1}_{\alpha i}\frac{\partial w^{i}}{\partial \bar{z}^{\beta}} d\bar{z}^{\beta}\otimes\frac{\partial}{\partial z^{\alpha}},
\end{equation}
where by abuse of notations, we write $w^i=f_t^*w^i=w^i\circ f_t$ for each $i=1,2,\cdots,n$. It can be checked that $\phi$ does not depend on the local coordinates and is a global vector form on $X$. Alternatively, $\phi(t):T^{0,1}\to T^{1,0}$ may be given by the following diagram:
\begin{equation}\label{phi}
\xymatrix{
                & T_{X,\mathbb{C}}  \ar[dl]_{p^{0,1}}      &&     T_{X,\mathbb{C}}  \ar[dr]^{p^{1,0}}        \\
-\phi(t): T_X^{0,1}  \ar[rr]^{(p^{0,1})^{-1}} && T_{X_t}^{0,1} \ar@{_{(}->}[lu]_{i} \ar@{^{(}->}[ru]^{j} \ar[rr]   &&   T_{X}^{1,0},      }
\end{equation}
where $i, j$ are inclusions and $p^{0,1}, p^{1,0}$ are projections.
But we note that the definition of $\phi(t)$ does depend on the choices of the smooth trivialization $F:X\times \Delta^k\to \mathcal{X}$ (whose existence is assured by the Ehresmann theorem) because we need to identify functions or vector fields on $X$ and on $X_t$. In what follows, we will always assume a trivialization $F$ is fixed.
\subsubsection{The exponential operator}
For any vector forms such as the Beltrami
differential $\phi(t)$, we can associate an operator $i_{\phi(t)}$ which is a derivation on the exterior algebra $A^\bullet(X)=\oplus_{p,q} A^{p,q}(X)$ of differential forms on $X$. For any $\varphi\in A^{p,q}(X)$,
\begin{equation}\label{i_phi definition}
i_{\phi(t)}\varphi:=\phi^{\alpha}_{\beta}d\bar{z}^{\beta}\wedge (\frac{\partial}{\partial z^{\alpha}}\lrcorner\varphi)\in A^{p-1,q+1}(X_0),
\end{equation}
where '$\lrcorner$' denotes the contraction operation (the operator $\frac{\partial}{\partial z^{\alpha}}\lrcorner \bullet$ is also called the interior derivative in the literature). We easily see that $(i_{\phi(t)})^{n+1} = 0$ so that
\begin{equation}\label{eq-exp}
e^{i_{\phi(t)}}:=\sum_{k=0}^{\infty} \frac{i_{\phi(t)}^k}{k!}~: A^\bullet(X) \to A^\bullet(X)
\end{equation}
is a well-defined operator which is called the \emph{exponential operator}. Since $e^{i_{\phi(t)}}e^{-i_{\phi(t)}}=e^{-i_{\phi(t)}}e^{i_{\phi(t)}}=e^0$ is the identity operator, $e^{-i_{\phi(t)}}$ is the inverse operator of $e^{i_{\phi(t)}}$. An important property of the exponential operator is that it preserves $(p,0)$-forms, that is, we have
\[
e^{i_{\phi(t)}}~: A^{p,0}(X)\longrightarrow A^{p,0}(X_t).
\]
Indeed, it follows from \eqref{Beltrami differential coord} that
\begin{equation}\label{dw in terms of dz}
d w^{\beta}=\frac{\partial w^{\beta}}{\partial z^{\alpha}}dz^{\alpha}+\frac{\partial w^{\beta}}{\partial \bar{z}^{\alpha}}d\bar{z}^{\alpha}
=\frac{\partial w^{\beta}}{\partial z^{\alpha}}(1+i_{\phi})dz^{\alpha}\in A^{1,0}(X_t),
\end{equation}
where $1$ denote the identity operator and
\[
e^{i_{\phi(t)}} (dz^{i_1}\wg\cdots\wg dz^{i_p})=(dz^{i_1}+i_{\phi(t)}dz^{i_1}) \wg \cdots \wg (dz^{i_p}+i_{\phi(t)}dz^{i_p}),
\]
which can be deduced from the fact that $e^{i_{\phi(t)}}$ is compatible with the wedge product (see e.g. \cite{Cle05}) and $e^{i_{\phi(t)}}dz^{i_k}=dz^{i_k}+i_{\phi(t)}dz^{i_k}$ for $k=1,\cdots,p$.
\subsubsection{The projection operator}
Let $z^1,\cdots, z^n$ and $w^1,\cdots, w^n$ be holomorphic coordinates on $X$ and $X_t$ respectively. Clearly, we have
\[
dz^{\alpha}=\frac{\partial z^{\alpha}}{\partial w^i} dw^i+ \frac{\partial z^{\alpha}}{\partial \bar{w}^i} d\bar{w}^i,
\]
and
\[
d\bar{z}^{\alpha}=\frac{\partial \bar{z}^{\alpha}}{\partial w^i} dw^i+ \frac{\partial \bar{z}^{\alpha}}{\partial \bar{w}^i} d\bar{w}^i.
\]
We can define an operator
\[
P~: A^{p,q}(X)\longrightarrow A^{p,q}(X_t),
\]
such that $P$ is compatible with the wedge product and
\[
Pf=f,\qquad P dz^{\alpha}=\frac{\partial z^{\alpha}}{\partial w^i} dw^i,\qquad P d\bar{z}^{\alpha}=\frac{\partial \bar{z}^{\alpha}}{\partial \bar{w}^i} d\bar{w}^i,
\]
where $f$ is a smooth function.
\begin{definition}\label{def-holo-tensor-bundle}
By a \emph{holomorphic tensor bundle} on a complex manifold $X$, we mean a holomorphic vector bundle formed by the tensor products or exterior products from the tangent bundle $T_{X}^{1,0}$ and its dual $\Omega_X=T_{X}^{1,0*}$.
\end{definition}
Let $E$ be a holomorphic tensor bundle on $X$ and $E_t$ the corresponding holomorphic tensor bundle on $X_t$. $P$ can be extended to an operator on $E$-valued $(0,q)$-forms,
\[
P~: A^{0,q}(X,E)\longrightarrow A^{0,q}(X_t,E_t),
\]
such that $P$ is compatible with the tensor product, wedge product and
\[
P\frac{\partial }{\partial z^{\alpha}}= \frac{\partial w^{i}}{\partial z^{\alpha}}\frac{\partial }{\partial w^{i}}.
\]
It can be checked easily that $P$ does not depend on the choices of local coordinates and is thus well-defined.
\subsubsection{The extension operator for $E$-valued $(0,q)$-forms}
\label{The extension operator for $E$-valued $(0,q)$-forms}
By combining the exponential operator $e^{i_{\phi(t)}}$ and the projection operator, we define a new operator
\begin{equation}\label{eq-ext-op}
\rho~:A^{0,q}(X, E)\to A^{0,q}(X_t, E_t),
\end{equation}
such that $\rho$ is compatible with the tensor product, wedge product and
\[
\rho\mid_{A^{0,q}(X)}=P,\qquad \rho\mid_{A^{0}(X,\Omega_X^p)}= e^{i_{\phi(t)}},\qquad \rho\mid_{A^{0}(X,T_X^{1,0})}=P.
\]
For example, if $E=T_{X}^{1,0}\wedge (T_{X}^{1,0} \otimes \Omega_X^p)$, we have
\begin{align*}
\rho:~~~&A^{0,q}(X, E)\longrightarrow A^{0,q}(X_t, E_t),\\
\sigma=\sigma_{\alpha\beta\gamma}\otimes\frac{\partial }{\partial z^{\alpha}}&\wedge (\frac{\partial }{\partial z^{\beta}}\otimes s^{\gamma}) \mapsto P\sigma_{\alpha\beta\gamma}\otimes \frac{\partial w^{i}}{\partial z^{\alpha}}\frac{\partial }{\partial w^{i}} \wedge (\frac{\partial w^{j}}{\partial z^{\beta}}\frac{\partial }{\partial w^{j}} \otimes e^{i_{\phi(t)}}s^{\gamma}),
\end{align*}
where $\sigma_{\alpha\beta\gamma}$ are local $(0,q)$-forms and $s^{\gamma}$ are local frames of $\Omega_X^p$. It can be checked that $\rho$ is a well-defined isomorphism. Note that for any open subset $U$ of $X$, $\rho:A^{0,q}(U, E)\to A^{0,q}(U, E_t)$ also make sense and hence $\rho$ is an isomorphism between the sheaf of germs of smooth $E$-valued $(0,q)$-forms on $X$ to the sheaf of germs of smooth $E$-valued $(0,q)$-forms on $X_t$.
\begin{remark}
In \cite{RZ15,RZ18}, Rao-Zhao defined an extension operator, denoted by $e^{i_{\phi(t)}\mid i_{\overline{\phi(t)}}}$, for $(p,q)$-forms as follows:
\begin{align*}
e^{i_{\phi(t)}\mid i_{\overline{\phi(t)}}}:~~~&A^{p,q}(X)\to A^{p,q}(X_t),\\
&\varphi=\varphi_{IJ}dz^I\wedge d\bar{z}^J\mapsto \varphi_{IJ}(e^{i_{\phi(t)}}dz^I)\wedge (e^{i_{\overline{\phi(t)}}}d\bar{z}^J),
\end{align*}
where $I = (i_1,\cdots, i_p)$, $J = (j_1,\cdots, j_q)$ are multi-index and $dz^I = dz^{i_1}\wedge\cdots\wedge dz^{i_p}$,
$d\bar{z}^J = d\bar{z}^{j_1} \wedge\cdots\wedge d\bar{z}^{j_q}$. An advantage of $e^{i_{\phi(t)}\mid i_{\overline{\phi(t)}}}$ is that it is a real operator~\cite[Lem.\,2.10]{RZ18} while $\rho\mid_{A^\bullet(X)}$ is not. But we will see that the extension equation for $\rho$ will be much simpler, see Theorem \ref{rho conjugated formula}.
\end{remark}
It follows from \eqref{dw in terms of dz} and the definition of $\rho$ that for any $\alpha,~i=1,\cdots,n$,
\begin{equation}\label{eq-inverse-rho-1}
\rho dz^{\alpha}=\left(\frac{\partial w}{\partial z}\right)^{-1}_{\alpha i}dw^{i},\quad~
\rho d\bar{z}^{\alpha}=\frac{\partial \bar{z}^{\alpha}}{\partial \bar{w}^i} d\bar{w}^i,
\end{equation}
and
\begin{equation}\label{eq-inverse-rho-2}
\rho^{-1} dw^{i}=\frac{\partial w^i}{\partial z^\alpha} dz^{\alpha},\quad~
\rho^{-1} d\bar{w}^i=\left(\frac{\partial \bar{z}}{\partial \bar{w}} \right)^{-1}_{i\alpha} d\bar{z}^{\alpha},\quad~
\rho^{-1}\frac{\partial }{\partial w^{i}}=\left(\frac{\partial w}{\partial z} \right)^{-1}_{\alpha i}\frac{\partial }{\partial z^{\alpha}}.
\end{equation}
For example, the last equality of \eqref{eq-inverse-rho-2} follows since
\[
\rho\frac{\partial }{\partial z^{\alpha}}=P\frac{\partial }{\partial z^{\alpha}}= \frac{\partial w^{i}}{\partial z^{\alpha}}\frac{\partial }{\partial w^{i}}.
\]
The following local computations will be useful:
\begin{lemma}\label{a local computation}
For any $t\in \Delta^k$, let $\bar{\partial}_t$ be the Dolbeault operator on $X_t$. Then we have
\begin{itemize}
\item[$(i)$] For $i=1,\cdots,n$,
\[
\frac{\partial }{\partial \bar{w}^{i}} = \frac{\partial \bar{z}^{\gamma}}{\partial \bar{w}^{i}}(\frac{\partial }{\partial \bar{z}^{\gamma}}- \phi^{\beta}_{\gamma} \frac{\partial }{\partial z^{\beta}}),
\]
and for any smooth function $f$ on $X$,
\[
\rho^{-1}\bar{\partial}_t\rho f= (\bar{\partial}- i_{\phi(t)}\partial)f;
\]
\item[$(ii)$] For any $\alpha=1,\cdots,n$,
\[
\rho^{-1}\bar{\partial}_t\rho d\bar{z}^{\alpha}=0,~\quad~\text{and}~\quad \rho^{-1}\bar{\partial}_t\rho dz^{\alpha}=\partial i_{\phi(t)}dz^{\alpha};
\]
\item[$(iii)$] For any $\alpha=1,\cdots,n$,
\[
\rho^{-1}\bar{\partial}_t\rho \frac{\partial }{\partial z^{\alpha}}=-[\phi(t),\frac{\partial }{\partial z^{\alpha}}].
\]
\end{itemize}
\end{lemma}

\begin{proof}For $(i)$, first we have
\begin{align*}
\frac{\partial }{\partial \bar{w}^{i}}&=\frac{\partial z^{\gamma}}{\partial \bar{w}^{i}}\frac{\partial }{\partial z^{\gamma}}
+ \frac{\partial \bar{z}^{\gamma}}{\partial \bar{w}^{i}}\frac{\partial }{\partial \bar{z}^{\gamma}}\\
&=\frac{\partial z^{\gamma}}{\partial \bar{w}^{i}}(\frac{\partial w^j}{\partial z^{\gamma}}\frac{\partial }{\partial w^{j}}+\frac{\partial \bar{w}^j}{\partial z^{\gamma}}\frac{\partial }{\partial \bar{w}^{j}})
+\frac{\partial \bar{z}^{\gamma}}{\partial \bar{w}^{i}}(\frac{\partial w^j}{\partial \bar{z}^{\gamma}}\frac{\partial }{\partial w^{j}}+\frac{\partial \bar{w}^j}{\partial \bar{z}^{\gamma}}\frac{\partial }{\partial \bar{w}^{j}}),\\
\end{align*}
which by comparing types implies
\begin{equation}\label{eq-phi-a-b}
-\left(\frac{\partial \bar{z}}{\partial \bar{w}}\right)^{-1}_{i\alpha}\frac{\partial z^{\beta}}{\partial \bar{w}^{i}} =\left(\frac{\partial w}{\partial z}\right)^{-1}_{\alpha i}\frac{\partial w^{i}}{\partial \bar{z}^{\beta}}=\phi^\alpha_\beta,\quad \forall\alpha,\beta=1,\cdots,n.
\end{equation}
Therefore,
\[
\frac{\partial }{\partial \bar{w}^{i}}=\frac{\partial z^{\gamma}}{\partial \bar{w}^{i}}\frac{\partial }{\partial z^{\gamma}}
+ \frac{\partial \bar{z}^{\gamma}}{\partial \bar{w}^{i}}\frac{\partial }{\partial \bar{z}^{\gamma}}
=\frac{\partial \bar{z}^{\gamma}}{\partial \bar{w}^{i}}(\frac{\partial }{\partial \bar{z}^{\gamma}} -\phi^{\beta}_{\gamma} \frac{\partial }{\partial z^{\beta}} ).
\]
Furthermore, it follows from \eqref{eq-inverse-rho-2} that
\[
\rho^{-1}\bar{\partial}_t\rho f=\rho^{-1}\frac{\partial f}{\partial \bar{w}^{i}}d\bar{w}^{i}=(\frac{\partial f }{\partial \bar{z}^{\alpha}} -\phi^{\beta}_{\alpha} \frac{\partial f}{\partial z^{\beta}} )d\bar{z}^\alpha
=(\bar{\partial}- i_{\phi(t)}\partial)f.
\]
For $(ii)$, it follows from \eqref{eq-inverse-rho-1} and \eqref{eq-inverse-rho-2} that
\[
\rho^{-1}\bar{\partial}_t\rho d\bar{z}^{\alpha}=\rho^{-1}\bar{\partial}_t\frac{\partial \bar{z}^{\alpha}}{\partial \bar{w}^i} d\bar{w}^i
=\rho^{-1}\bar{\partial}_t\bar{\partial}_t (\bar{z}^{\alpha})=0.
\]
and
\begin{align*}
 \rho^{-1}\bar{\partial}_t\rho dz^{\alpha}
&=\rho^{-1}\bar{\partial}_t(1+i_{\phi})dz^{\alpha}\\
&=\rho^{-1}\bar{\partial}_t\left(\frac{\partial w}{\partial z}\right)^{-1}_{\alpha i}dw^{i}\\
&=[\frac{\partial }{\partial \bar{w}^{j}}\left(\frac{\partial w}{\partial z}\right)^{-1}_{\alpha i}]~\rho^{-1}d\bar{w}^{j}\wedge \rho^{-1}dw^{i}\\
&=[\frac{\partial }{\partial \bar{w}^{j}}\left(\frac{\partial w}{\partial z}\right)^{-1}_{\alpha i}]\left(\frac{\partial \bar{z}}{\partial \bar{w}}\right)^{-1}_{j \nu}\frac{\partial w^{i}}{\partial z^{\tau}} d\bar{z}^{\nu}\wedge dz^{\tau}\\
&=-[\frac{\partial }{\partial \bar{w}^{j}} \frac{\partial w^{i}}{\partial z^{\tau}}) ]\left(\frac{\partial \bar{z}}{\partial \bar{w}}\right)^{-1}_{j \nu} \left(\frac{\partial w}{\partial z}\right)^{-1}_{\alpha i} d\bar{z}^{\nu}\wedge dz^{\tau}\\
&=-[\frac{\partial \bar{z}^{\gamma}}{\partial \bar{w}^{j}} (\frac{\partial }{\partial \bar{z}^{\gamma}} \frac{\partial w^{i}}{\partial z^{\tau}} - \phi^{\beta}_{\gamma} \frac{\partial }{\partial z^{\beta}} \frac{\partial w^{i}}{\partial z^{\tau}}) ]\left(\frac{\partial \bar{z}}{\partial \bar{w}}\right)^{-1}_{j \nu} \left(\frac{\partial w}{\partial z}\right)^{-1}_{\alpha i} d\bar{z}^{\nu}\wedge dz^{\tau}\\
&=-[\frac{\partial }{\partial z^{\tau}} (\frac{\partial w^{i}}{\partial \bar{z}^{\gamma}} - \phi^{\beta}_{\gamma} \frac{\partial w^{i}}{\partial z^{\beta}}) +\frac{\partial \phi^{\beta}_{\gamma}}{\partial z^{\tau}}\frac{\partial w^{i}}{\partial z^{\beta}} ] \left(\frac{\partial w}{\partial z}\right)^{-1}_{\alpha i} d\bar{z}^{\gamma}\wedge dz^{\tau}\\
&=\partial \phi^{\alpha}_{\gamma}\wedge d\bar{z}^{\gamma}=\partial i_{\phi(t)}dz^{\alpha},\\
\end{align*}
where we have used the fact that
\[
\frac{\partial w^{i}}{\partial \bar{z}^{\gamma}} - \phi^{\beta}_{\gamma} \frac{\partial w^{i}}{\partial z^{\beta}}=0,\quad~\text{for}~i,\gamma=1,\cdots,n.
\]
For $(iii)$,
\begin{align*}
 \rho^{-1}\bar{\partial}_t\rho \frac{\partial }{\partial z^{\alpha}}
&=\rho^{-1}\bar{\partial}_t \frac{\partial w^{i}}{\partial z^{\alpha}} \frac{\partial }{\partial w^{i}}\\
&=\rho^{-1} \frac{\partial }{\partial \bar{w}^{j}} \frac{\partial w^{i}}{\partial z^{\alpha}} d\bar{w}^{j}\otimes \frac{\partial }{\partial w^{i}}\\
&=\frac{\partial }{\partial \bar{w}^{j}} \frac{\partial w^{i}}{\partial z^{\alpha}} \left(\frac{\partial \bar{z}}{\partial \bar{w}} \right)^{-1}_{j\gamma} d\bar{z}^{\gamma} \otimes \left(\frac{\partial w}{\partial z}\right)^{-1}_{\alpha i} \frac{\partial }{\partial z^{\alpha}} \\
&= (\frac{\partial }{\partial \bar{z}^{\beta}} \frac{\partial w^{i}}{\partial z^{\alpha}} - \phi^{\gamma}_{\beta} \frac{\partial }{\partial z^{\gamma}} \frac{\partial w^{i}}{\partial z^{\alpha}} ) d\bar{z}^{\beta} \otimes \left(\frac{\partial w}{\partial z}\right)^{-1}_{\alpha i} \frac{\partial }{\partial z^{\alpha}} \\
&=[ \frac{\partial }{\partial z^{\alpha}} ( \frac{\partial w^{i}}{\partial \bar{z}^{\beta}} - \phi^{\gamma}_{\beta} \frac{\partial w^{i}}{\partial z^{\gamma}} ) + \frac{\partial \phi^{\gamma}_{\beta}}{\partial z^{\alpha}} \frac{\partial w^{i}}{\partial z^{\gamma}}] d\bar{z}^{\beta} \otimes  \left(\frac{\partial w}{\partial z}\right)^{-1}_{\alpha i} \frac{\partial }{\partial z^{\alpha}} \\
&=\frac{\partial \phi^{\gamma}_{\beta}}{\partial z^{\alpha}} d\bar{z}^{\beta} \otimes \frac{\partial }{\partial z^{\gamma}}=-[\phi(t),\frac{\partial }{\partial z^{\alpha}}].
\end{align*}
\end{proof}
\subsection{The extension equation}
\label{The extension equation}
First, we recall some basic facts about derivations. See \cite{KMS93,FN56} and \cite{Xia19deri} for more details.
\subsubsection{Derivations and bracket operations on a real manifold }
Let $M$ be a smooth manifold of dimension $n$, $A(M)=\bigoplus_{k=0}^{k=n} A^k(M)$ be its exterior algebra of differential forms. A (graded) \emph{derivation} $D$ on $A(M)$ is a $\mathbb{R}$-linear map $D:A(M)\to A(M)$ with $D(A^l(M))\subseteq A^{l+k}(M)$ and $D(\xi\wedge\eta)=(D\xi)\wedge\eta+(-1)^{kl}\xi\wedge(D\eta)$ for $\xi\in A^l(M)$. The integer $k$ is called the degree of $D$. We denote by $D^k(M)$ the space of all derivations of degree $k$ on $A(M)$. For any $D_1\in D^{k_1}(M)$ and $D_2\in D^{k_2}(M)$, the graded commutator is defined by $[D_1,D_2]:=D_1D_2-(-1)^{k_1k_2}D_2D_1$. With respect to this bracket operation, the space of all derivations $D(M)=\bigoplus_kD^k(M)$ becomes a graded Lie algebra.

For a vector valued $(k+1)$-form $K\in A^{k+1}(M,TM)$, we can associate a derivation $i_K$ of degree $k$ by setting $i_K\varphi:=\xi\wedge(i_X\varphi)$, if $K=\xi\otimes X$ for a $(k+1)$-form $\xi$ and a vector field $X$, where $\varphi\in A(M)$. Let $d$ be the exterior derivative on $M$, then the Lie derivative is defined as $\mathcal{L}_K:=[i_K,d]=i_Kd-(-1)^{k}di_K$, where $K\in A^{k+1}(M,TM)$. According to a result due to Fr\"olicher-Nijenhuis~\cite[Prop.\,4.7]{FN56}, for any $D\in D^k(M)$, there exist unique $K\in A^{k}(M,TM)$ and $L\in A^{k+1}(M,TM)$ such that \begin{equation}\label{D=lK+iL}
D=\mathcal{L}_K+i_L,
\end{equation}
and $L=0$ if and only if $[D,d]=0$, $K=0$ if and only if $D$ is algebraic.

A derivation $D$ on $A(M)$ is called \emph{algebraic} if $Df=0, \forall f\in A^0(M)$. By \eqref{D=lK+iL}, every algebraic derivation of degree $k$ on $A(M)$ is of the form $i_L$ for some unique $L\in A^{k+1}(M,TM)$.

For any $K\in A^{k}(M,TM)$, $L\in A^{l}(M,TM)$, $[\mathcal{L}_K,\mathcal{L}_L]$ is a derivation of degree $k+l$ such that $[[\mathcal{L}_K,\mathcal{L}_L],d]=0$, hence by \eqref{D=lK+iL},
\begin{equation}\label{ll commutator}
[\mathcal{L}_K,\mathcal{L}_L]=\mathcal{L}_{[K,L]}
\end{equation}
for some unique $[K,L]\in A^{k+l}(M,TM)$. This operation $[\cdot,\cdot]$ is called the \emph{Fr\"{o}licher-Nijenhuis bracket}.
\subsubsection{Lie derivatives on complex manifolds}
Let $X$ be a complex manifold, we have the following
\begin{lemma}\cite{Xia19deri}\label{alg K}
For any vector $k$-form $K$ on $X$, set $\mathcal{L}_{K}^{1,0}:=[i_K,\partial]=i_K\partial-(-1)^{k-1}\partial i_K$, and $\mathcal{L}_{K}^{0,1}:=[i_K,\bar{\partial}]=i_K\bar{\partial}-(-1)^{k-1}\bar{\partial} i_K$, we have
\begin{itemize}
  \item[(1)]if $K\in A^{0,k}(X,T^{1,0})$, then $\mathcal{L}_{K}^{0,1}= (-1)^{k}i_{\bar{\partial}K}$;
  \item[(2)] if $K\in A^{k,0}(X,T^{0,1})$, then $\mathcal{L}_{K}^{1,0}= (-1)^{k}i_{\partial K}$;
  \item[(3)] if $K\in A^{k,0}(X,T^{1,0})$, then $\mathcal{L}_{K}^{0,1}= (-1)^{k}i_{\bar{\partial}K}$;
  \item[(4)] if $K\in A^{0,k}(X,T^{0,1})$, then $\mathcal{L}_{K}^{1,0}= (-1)^{k}i_{\partial K}$.
\end{itemize}
\end{lemma}
\begin{lemma}\label{used in thm}
 Let $\varphi\in A^{0,1}(X,T^{1,0})$ and $\psi\in A^{0,1}(X,T^{1,0})$, then we have
\begin{itemize}
  \item[(1)] $[\mathcal{L}_{\varphi}^{1,0},\mathcal{L}_{\psi}^{1,0}]=\mathcal{L}_{[\varphi,\psi]}^{1,0}$;
  \item[(2)] $[\bar{\partial}, \mathcal{L}_{\varphi}^{1,0}]=\mathcal{L}_{\bar{\partial}\varphi}^{1,0}$.
\end{itemize}
\end{lemma}

\begin{proof}By \eqref{ll commutator} and Lemma \ref{alg K}, we have
\[
[\mathcal{L}_{\varphi},\mathcal{L}_{\psi}]=\mathcal{L}_{[\varphi,\psi]}=\mathcal{L}_{[\varphi,\psi]}^{1,0}+i_{\bar{\partial}[\varphi, \psi]},
\]
on the other hand,
\[
[\mathcal{L}_{\varphi},\mathcal{L}_{\psi}]=[\mathcal{L}_{\varphi}^{1,0}-i_{\bar{\partial}\varphi}, \mathcal{L}_{\psi}^{1,0}-i_{\bar{\partial}\psi}]
=[\mathcal{L}_{\varphi}^{1,0},\mathcal{L}_{\psi}^{1,0}] - [\mathcal{L}_{\varphi}^{1,0}, i_{\bar{\partial}\psi}] - [i_{\bar{\partial}\varphi}, \mathcal{L}_{\psi}^{1,0}] + [i_{\bar{\partial}\varphi}, i_{\bar{\partial}\psi}].
\]
By comparing the bidegree, $(1)$ follows. For $(2)$, by applying the Jacobi identity and Lemma \ref{alg K}, we have
\[[\bar{\partial}, \mathcal{L}_{\varphi}^{1,0}]=[\bar{\partial}, [i_{\varphi}, \partial]]=[[\bar{\partial}, i_{\varphi}], \partial]+ [i_{\varphi}, [\bar{\partial}, \partial]]=[-\mathcal{L}_{\varphi}^{0,1}, \partial]=[i_{\bar{\partial}\varphi}, \partial]=\mathcal{L}_{\bar{\partial}\varphi}^{1,0}.
\]
\end{proof}
Formulas in the above two lemmas are also called the holomorphic Cartan homotopy formulas, see \cite{BM18,FM09,FM06,Cle05} for more information.
\subsubsection{The case of holomorphic exterior bundle}
Let $\pi : \mathcal{X} \to \Delta^k$ be a complex analytic family with Beltrami differentials $\phi(t)$, we have
\begin{theorem}\label{rho conjugated formula} Let $d=\partial+ \bar{\partial}= \partial_t+ \bar{\partial}_t$, where $\bar{\partial}$ and $\bar{\partial}_t$ are the Dolbeault operators on
$X$ and $X_t$, respectively. Then we have
\begin{equation}\label{rho conjugated p-bar formula}
\rho^{-1}\bar{\partial}_t\rho = \bar{\partial}-\mathcal{L}_{\phi(t)}^{1,0}~\quad\text{on}\quad A^\bullet(X).
\end{equation}
\end{theorem}
\begin{proof}Since this is an equality of derivations on $A^\bullet(X)$, it is enough to prove \eqref{rho conjugated p-bar formula} for functions and $1$-forms which is a direct consequence of the first two assertions of Lemma \ref{a local computation}.
\end{proof}

\subsubsection{The case of tensor products of holomorphic exterior bundle}
Before consider general holomorphic tensor bundles, we first examine carefully the case of tensor products of holomorphic exterior bundle which will help us understand the general situation. Let $\pi : \mathcal{X} \to \Delta^k$ be an analytic family of compact complex manifolds of dimension $n$ over the unit polydisc $\Delta^k$ in $\mathbb{C}^k$ with Beltrami differentials $\phi(t)$. Set $E_t:=\Omega_{X_t}^{p_1}\otimes\cdots\otimes\Omega_{X_t}^{p_m}$, $\forall t \in \Delta^k$, where $\Omega_{X_t}^{p}:=\wedge^{p}T^{*1,0}_{X_t}$. Given a smooth global section $\sigma$ of the bundle $\wedge^{q}T^{*0,1}_{X_0}\otimes E_0$, i.e. $\sigma\in A^{0,q}(X_0,E_0)$, locally we can write $\sigma=\varphi_{\alpha_1\cdots\alpha_m}\otimes e_{\alpha_1}\otimes\cdots\otimes e_{\alpha_m}$, where $\left\{e_{\alpha_i}\right\}$ are local smooth frame of $\wedge^{p_i}T^{*1,0}_{X_0}=\Omega_{X_0}^{p_i}$, $i=1, \cdots, m$, and $\varphi_{\alpha_1\cdots\alpha_m}$ are $(0,q)$-forms. In this case, the extension operator $\rho$ is given by:
\begin{align*}
\rho:~~~&A^{0,q}(X,E)\longrightarrow A^{0,q}(X_t,E_t),\\
\sigma=\varphi_{\alpha_1\cdots\alpha_m}\otimes& e_{\alpha_1}\otimes\cdots\otimes e_{\alpha_m} \longmapsto P\varphi_{\alpha_1\cdots\alpha_m}\otimes e^{i_{\phi(t)}} e_{\alpha_1}\otimes\cdots\otimes e^{i_{\phi(t)}} e_{\alpha_m}.
\end{align*}

Let $E$ be a holomorphic vector bundle on $X=X_0$, the Dolbeault operator $\bar{\partial}$ can be naturally defined on $A^{0,q}(X,E)$ such that
\[\bar{\partial}: A^{0,q}(X,E)\to A^{0,q+1}(X,E).
\]
In fact, given $s\in A^{0,q}(X,E)$, we can write $s=s_{\alpha}e_{\alpha}$ in terms of a holomorphic frame $\left\{e_{\alpha}\right\}$, and define \begin{equation}\label{p-bar def}
\bar{\partial}s :=\bar{\partial}s_{\alpha}\otimes e_{\alpha}.
\end{equation}It can be checked that this definition does not depend on the holomorphic frame chosen~\cite{GH94}. If we write $s=s'_{\alpha}e'_{\alpha}$ where $\left\{e'_{\alpha}\right\}$ is merely a smooth frame, then we would have
\[\bar{\partial}s=\bar{\partial}s'_{\alpha}\otimes e'_{\alpha} + (-1)^q s'_{\alpha}\otimes \bar{\partial}e'_{\alpha},
\]
where $\bar{\partial}e'_{\alpha}$ is defined by \eqref{p-bar def}. Motivated by this observation, we now extend the definition of the operators in previous sections to holomorphic tensor bundles. Consider the Lie derivative $\mathcal{L}_{K}^{1,0}$, where $K\in A^{0,1}(X,T_{X}^{1,0})$. Given $\sigma\in A^{0,q}(X,E)$, locally we write $\sigma=\varphi_{\alpha_1\cdots\alpha_m}\otimes e_{\alpha_1}\otimes\cdots\otimes e_{\alpha_m}$ as above and define $\tilde{\mathcal{L}}_{K}^{1,0}: A^{0,q}(X,E)\to A^{0,q+1}(X,E)$ by
\begin{equation}\label{Lie derivative extended}
\tilde{\mathcal{L}}_{K}^{1,0}\sigma=\mathcal{L}_{K}^{1,0}\varphi_{\alpha_1\cdots\alpha_m}\otimes e_{\alpha_1}\otimes\cdots\otimes e_{\alpha_m}+
(-1)^q\sum_{i=1}^{m}\varphi_{\alpha_1\cdots\alpha_m}\otimes e_{\alpha_1}\otimes\cdots\otimes \mathcal{L}_{K}^{1,0}e_{\alpha_i} \otimes\cdots\otimes e_{\alpha_m}.
\end{equation}

\begin{lemma} \eqref{Lie derivative extended} give rise to a well-defined global operator
\[
\tilde{\mathcal{L}}_{K}^{1,0}: A^{0,q}(X,E)\to A^{0,q+1}(X,E).
\]
\end{lemma}

\begin{proof}We will only prove this in the case $E=K^2_{X}$, where $K_{X}$ is the canonical line bundle on $X$, because the general case is essentially the same. Indeed, let $\sigma\in A^{0,q}(X,K^2_{X})$, locally we write $\sigma=\varphi\otimes e\otimes e$ in terms of the local frame $e$ of $K_{X}$. Let $e'$ be another local frame and $e=ge'$ where $g$ is a nonvanishing smooth function, then if $\sigma=\varphi'\otimes e'\otimes e'$, we have $\varphi'=g^2\varphi$ and
\begin{align*}
&\tilde{\mathcal{L}}_{K}^{1,0}(\varphi\otimes e\otimes e)\\
=&\mathcal{L}_{K}^{1,0}\varphi\otimes e\otimes e + (-1)^q(\varphi\otimes \mathcal{L}_{K}^{1,0}e\otimes e + \varphi\otimes e\otimes \mathcal{L}_{K}^{1,0}e)\\
=&-2\frac{1}{g}\mathcal{L}_{K}^{1,0}g\wedge\varphi'\otimes e'\otimes e' + \mathcal{L}_{K}^{1,0}\varphi'\otimes e'\otimes e' +\\
&(-1)^q[\frac{1}{g}\varphi'\otimes (\mathcal{L}_{K}^{1,0}g\otimes e'+g\mathcal{L}_{K}^{1,0}e')\otimes e'+ \frac{1}{g}\varphi'\otimes e' \otimes (\mathcal{L}_{K}^{1,0}g\otimes e'+g\mathcal{L}_{K}^{1,0}e')]\\
=&-2\frac{1}{g}\mathcal{L}_{K}^{1,0}g\wedge\varphi'\otimes e'\otimes e' +  \mathcal{L}_{K}^{1,0}\varphi'\otimes e'\otimes e' +\\
 &(-1)^q[\frac{2}{g}\varphi'\wedge \mathcal{L}_{K}^{1,0}g\otimes e'\otimes e'+ \varphi' \otimes\mathcal{L}_{K}^{1,0}e'\otimes e' + \varphi' \otimes e'\otimes \mathcal{L}_{K}^{1,0} e']\\
=& \mathcal{L}_{K}^{1,0}\varphi'\otimes e'\otimes e'+ (-1)^q(\varphi' \otimes\mathcal{L}_{K}^{1,0}e'\otimes e' + \varphi' \otimes e'\otimes \mathcal{L}_{K}^{1,0} e').
\end{align*}
This shows that the definition of $\tilde{\mathcal{L}}_{K}^{1,0}$ does not depend on the local frame chosen.
\end{proof}

\begin{theorem}\label{ext formula for tensor product of exterior bundle}
 Let $\bar{\partial}$ and $\bar{\partial}_t$ be the Dolbeault operator on $E$ and $E_t$, respectively. Then we have
\begin{equation}
\rho^{-1}\bar{\partial}_t\rho = \bar{\partial}-\tilde{\mathcal{L}}_{\phi(t)}^{1,0}~\quad\text{on}\quad A^{0,\bullet}(X,E).
\end{equation}
In particular, for $\sigma\in A^{0,q}(X,E)$, $\rho\sigma \in A^{0,q}(X_t,E_t)$ is $\bar{\partial}_t$-closed if and only if
\begin{equation}
(\bar{\partial}-\tilde{\mathcal{L}}_{\phi(t)}^{1,0}) \sigma=0.
\end{equation}
\end{theorem}

\begin{proof}Given $\sigma\in A^{0,q}(X,E)$, locally we write $\sigma=\varphi_{\alpha_1\cdots\alpha_m}\otimes e_{\alpha_1}\otimes\cdots\otimes e_{\alpha_m}$, then it follows from Theorem \ref{rho conjugated formula} that
\begin{align*}
\rho^{-1}\bar{\partial}_t\rho \sigma
=&(\rho^{-1}\bar{\partial}_t\rho)\varphi_{\alpha_1\cdots\alpha_m}\otimes e_{\alpha_1}\otimes\cdots\otimes e_{\alpha_m}+\\
&(-1)^q\sum_{i=1}^{m}\varphi_{\alpha_1\cdots\alpha_m}\otimes e_{\alpha_1}\otimes\cdots\otimes (\rho^{-1}\bar{\partial}_t\rho) e_{\alpha_i} \otimes\cdots\otimes e_{\alpha_m}\\
=&(\bar{\partial}-\mathcal{L}_{\phi(t)}^{1,0})\varphi_{\alpha_1\cdots\alpha_m}\otimes e_{\alpha_1}\otimes\cdots\otimes e_{\alpha_m}+\\
&(-1)^q\sum_{i=1}^{m}\varphi_{\alpha_1\cdots\alpha_m}\otimes e_{\alpha_1}\otimes\cdots\otimes (\bar{\partial}-\mathcal{L}_{\phi(t)}^{1,0}) e_{\alpha_i} \otimes\cdots\otimes e_{\alpha_m}\\
=&(\bar{\partial}-\tilde{\mathcal{L}}_{\phi(t)}^{1,0}) \sigma.
\end{align*}
\end{proof}
\subsubsection{The extension equation for $E$-valued $(0,q)$-forms}
In Subsection \ref{The extension operator for $E$-valued $(0,q)$-forms}, we have defined the extension operator for $E$-valued $(0,q)$-forms. We are now ready to derive the extension equation in this general situation.

Let $\pi : \mathcal{X} \to \Delta^k$ be a complex analytic family with Beltrami differentials $\phi(t)$. Let $E$ be a holomorphic tensor bundle on $X$ and $E_t$ the corresponding holomorphic tensor bundle on $X_t$. Clearly, we have the following commutative diagram
\[
\xymatrix{
  A^{0,q}(X, E) \ar[d]_{\rho^{-1}\bar{\partial}_t\rho} \ar[r]^{\rho} & A^{0,q}(X_t, E_t) \ar[d]^{\bar{\partial}_t} \\
  A^{0,q+1}(X, E) \ar[r]^{\rho} & A^{0,q+1}(X_t, E_t),   }
\]
where $E_t$ is the corresponding holomorphic tensor bundle on $X_t$.

Let $\phi\in A^{0,k}(X,T_{X}^{1,0})$. It follows directly form the definition of the Fr\"{o}licher-Nijenhuis bracket that for any $\varphi\in A^{0,p}(X)$ and $\sigma\in A^{0,\bullet}(X,T_{X}^{1,0})$, we have
\begin{equation}\label{2.24}
[\phi, \varphi\wedge\sigma ]=\mathcal{L}_{\phi}^{1,0}\varphi + (-1)^{pk}\varphi\wedge [\phi, \sigma ].
\end{equation}
Hence, by letting the operator $\langle\phi|$ satisfy the Leibniz's rule (as in \eqref{Lie derivative extended}) and set
\begin{equation}\label{eq-langle-phi}
\langle\phi|:=\left\{
\begin{array}{ll}
[\phi, \cdot],           &~~~~~~~~~~~~\text{on}~~~~~~~~~~A^{0,\bullet}(X, T_{X}^{1,0}),  \\
\mathcal{L}_{\phi}^{1,0},&~~~~~~~~~~~~\text{on}~~~~~~~~~~A^{0,\bullet}(X, \Omega_X), \\
\end{array} \right.
\end{equation}
we get a well-defined operator
\[
\langle\phi|~:A^{0,\bullet}(X, E)  \longrightarrow  A^{0,\bullet+k}(X, E)
\]
and a natural paring
\begin{align*}
\langle\cdot|\cdot\rangle : A^{0,k}(X, T_{X}^{1,0})\times A^{0,q}(X, E)  \longrightarrow & A^{0,q+k}(X, E),\\
(\phi , \sigma)\longmapsto & \langle\phi|\sigma\rangle:=\langle\phi|\sigma.
\end{align*}
We remark that by Leibniz's rule, to show the paring $\langle\cdot|\cdot\rangle$ satisfy some property it is enough to prove this when $E=T_{X}^{1,0}$ or $E=\Omega_{X}$. For example, let $\phi\in A^{0,k}(X, T_{X}^{1,0})$ and $\psi\in A^{0,l}(X, T_{X}^{1,0})$, it follows from Lemma \ref{used in thm}, (1) that
\begin{equation}\label{phi-psi-neat}
[\langle\phi|, \langle\psi|]=\langle\phi| \langle\psi| -(-1)^{kl}\langle\psi| \langle\phi| =\langle [\phi, \psi]| ,
\end{equation}
or
\begin{equation}\label{phi-psi}
\langle\phi| \langle\psi| \sigma\rangle\rangle =  \langle [\phi, \psi]|\sigma\rangle + (-1)^{kl}\langle\psi| \langle\phi|\sigma\rangle\rangle,~~~~~~~~~~\forall \sigma\in A^{0,\bullet}(X, E).
\end{equation}
In particular, for $\phi\in A^{0,1}(X, T_{X}^{1,0})$, we have
\begin{equation}\label{phi-phi}
\langle\phi| \langle\phi| \sigma\rangle\rangle =  \frac{1}{2}\langle [\phi, \phi]|\sigma\rangle ,~~~~~~~~~~\forall \sigma\in A^{0,\bullet}(X, E).
\end{equation}

\begin{theorem}\label{ext formula for general holo tensor bundle}
Let $\bar{\partial}$ and $\bar{\partial}_t$ be the Dolbeault operator on $X$ and $X_t$, respectively. Then we have
\begin{equation}
\rho^{-1}\bar{\partial}_t\rho = ~\bar{\partial}- \langle\phi(t) | ,~~~~~~~~\text{on}~~~~~~A^{0,\bullet}(X, E).
\end{equation}
In particular, for $\sigma\in A^{0,q}(X, E)$, $\rho \sigma \in A^{0,q}(X_t, E_t)$ is $\bar{\partial}_t$-closed if and only if
\begin{equation}
\bar{\partial}\sigma - \langle \phi(t) | \sigma \rangle =0.
\end{equation}
\end{theorem}
\begin{proof}This follows immediately from the same proof as Theorem \ref{ext formula for tensor product of exterior bundle} and the third assertion of Lemma \ref{a local computation}.
\end{proof}

\begin{remark}
In the case of $E=T^{1,0}$, Theorem \ref{ext formula for general holo tensor bundle} says that
\begin{equation}
\rho^{-1}\bar{\partial}_t\rho = \bar{\partial}- [\phi(t) , \cdot],~~~~~~~~\text{on}~~~~~~A^{0,\bullet}(X, T_{X}^{1,0}).
\end{equation}
This was first proved in~\cite{Ham77}, see also~\cite{Hua95,Zha14}.
\end{remark}
\subsection{DGLA-module and Extension isomorphism }

\subsubsection{DGLA-module}

Recall that~\cite{Ma05} a \emph{differential graded Lie algebra} (DGLA for short) is the data consists of a $\mathbb{Z}$-graded vector space $L=\oplus_{i\in \mathbb{Z}}L^i$ together with a bilinear bracket $[\cdot,\cdot]: L\times L \to L$, and a linear map $d: L \to L$ such that
\begin{enumerate}
  \item $[L^i, L^j]\subset L^{i+j}$ and $[a, b]= -(-1)^{ij}[b, a]$, if $a\in L^i$, $b\in L^j$;
  \item graded Jacobi identity: $[a,[b,c]]=[[a,b],c]+(-1)^{ij}[b,[a,c]]$, $\forall a\in L^i$, $b\in L^j$;
  \item $dL^i\subset L^{i+1}, d^2=d\circ d=0$, and $d[a, b]=[da, b]+ (-1)^{i}[a, db], \forall a\in L^i$.
\end{enumerate}
We will denote a DGLA by $(L, d, [\cdot,\cdot])$. For example, in classical deformation theory we have the Kodaira-Spencer DGLA $(A^{0,\bullet}(X, T_{X}^{1,0}), \bar{\partial}, [\cdot,\cdot])$. If we define $\bar{\partial}_{\phi(t)}:= \bar{\partial}- [\phi(t) , \cdot]$, then
\[
\bar{\partial}_{\phi(t)}:A^{0,\bullet}(X, T_{X}^{1,0})\to A^{0,\bullet+1}(X, T_{X}^{1,0}),
\]
and for each $t\in \Delta^k$, $(A^{0,\bullet}(X, T_{X}^{1,0}), \bar{\partial}_{\phi(t)}, [\cdot,\cdot])$ is a DGLA (see Proposition \ref{DGLA-mod-discription}), which will be called the \emph{deformed Kodaira-Spencer DGLA}.
\begin{definition}\cite{Gri64}\label{def-DGLA-module}
Let $(L, d_L, [\cdot,\cdot])$ be a DGLA. By a \emph{DGLA-module over} $(L, d_L, [\cdot,\cdot])$ we mean a complex $(V, d_V)$ together with a bilinear paring
\[
\langle\cdot, \cdot\rangle : L^i\times V^j \longrightarrow V^{i+j}
\]
such that
\[
d_V \langle a, x\rangle = \langle d_La, x\rangle + (-1)^{i}\langle a, d_Vx\rangle,~~~~~~~~~~\forall a\in L^i, x\in V.
\]
A DGLA-module will be denoted by $(V, d_V, \langle\cdot, \cdot\rangle)$.
\end{definition}

Let $X$ be a complex manifold and $X_t$ a small deformation (of $X$) whose complex structure is represented by a Beltrami differential $\phi(t)\in A^{0,1}(X, T_{X}^{1,0})$. For any holomorphic tensor bundle $E$ on $X$, consider the following operator
\[
\bar{\partial}_{\phi(t)}:= \bar{\partial}-\langle\phi(t)| : A^{0,\bullet}(X, E)\longrightarrow A^{0,\bullet+1}(X, E).
\]
We have the following observation:

\begin{proposition}\label{DGLA-mod-discription}
$\bar{\partial}_{\phi(t)}^2 = \bar{\partial}_{\phi(t)}\circ \bar{\partial}_{\phi(t)} =0$ and $\forall \varphi\in A^{0,k}(X, T_{X}^{1,0}), \sigma\in A^{0,\bullet}(X, E)$,
\begin{equation}\label{DGLA-mod-compatible}
\bar{\partial}_{\phi(t)} \langle \varphi| \sigma\rangle = \langle \bar{\partial}_{\phi(t)} \varphi| \sigma\rangle + (-1)^{k}\langle \varphi| \bar{\partial}_{\phi(t)}\sigma\rangle,
~~~\text{or}~~~[\bar{\partial}_{\phi(t)}, \langle \varphi| ]=\langle \bar{\partial}_{\phi(t)}\varphi|.
\end{equation}
In other words, $(A^{0,\bullet}(X, E), \bar{\partial}_{\phi(t)} )$ is a complex and $(A^{0,\bullet}(X, E), \bar{\partial}_{\phi(t)}, \langle\cdot|\cdot\rangle )$ is a DGLA-module over the deformed Kodaira-Spencer DGLA $(A^{0,\bullet}(X, T_{X}^{1,0}), \bar{\partial}_{\phi(t)}, [\cdot,\cdot])$.
\end{proposition}

\begin{proof}For \eqref{DGLA-mod-compatible}, by Leibniz's rule we only need to prove it for $E=T_{X}^{1,0}$ or $E=\Omega_{X}$. Since $\bar{\partial}_{\phi(t)}=\bar{\partial}-\langle\phi(t)|$ , in view of \eqref{phi-psi-neat}, it suffices to show $[\bar{\partial}, \langle \varphi| ]=\langle \bar{\partial}\varphi|$. For $[\bar{\partial}, \langle \varphi| ]=\langle \bar{\partial}\varphi|$, the case $E=T_{X}^{1,0}$ follows from the standard fact that
\[
\bar{\partial}[\varphi,\bullet]=[\bar{\partial}\varphi,\bullet]+(-1)^k[\varphi,\bar{\partial}\bullet],
\]
while the case $E=\Omega_{X}$ follows from Lemma \ref{used in thm}, (2). Given $\sigma\in A^{0,\bullet}(X, E)$, we have
\begin{align*}
\bar{\partial}_{\phi(t)}^2\sigma
&=\bar{\partial}_{\phi(t)} (\bar{\partial}\sigma - \langle \phi(t)| \sigma\rangle)\\
&=-\bar{\partial}\langle \phi(t)| \sigma\rangle - \langle \phi(t)| \bar{\partial}\sigma\rangle + \langle \phi(t)| \langle \phi(t)| \sigma\rangle\rangle \\
&=-\langle \bar{\partial}\phi(t)| \sigma\rangle + \frac{1}{2}\langle [\phi(t), \phi(t)]|\sigma\rangle\\
&=0,
\end{align*}
where we have used \eqref{DGLA-mod-compatible}, \eqref{phi-phi} and the Maurer-Cartan equation $\bar{\partial}\phi(t) = \frac{1}{2}[\phi(t) ,\phi(t) ]$.
\end{proof}
The complex $(A^{0,\bullet}(X, E), \bar{\partial}_{\phi(t)} )$ will be called the \emph{deformed Dolbeault complex}. the The cohomology of the complex $(A^{0,\bullet}(X, E), \bar{\partial}_{\phi(t)} )$ is denoted by
\[
H^{0,q}_{\bar{\partial}_{\phi(t)}}(X,E) := \frac{\ker \bar{\partial}_{\phi(t)} \cap A^{0,q}(X,E)}{\Image \bar{\partial}_{\phi(t)} \cap A^{0,q}(X,E)},~\qquad \forall q\geq 0 ,
\]
and we call $H^{0,q}_{\bar{\partial}_{\phi(t)}}(X,E)$ the \emph{deformed Dolbeault cohomology group}.
\subsubsection{Extension isomorphism}\label{sec-Extension isomorphism}
Now it is a direct consequence of Theorem \ref{ext formula for general holo tensor bundle} that the following commutative diagram is commutative:
\[
\xymatrix{
  A^{0,\bullet}(X, E) \ar[d]_{\bar{\partial}_{\phi(t)} } \ar[r]^{\rho} & A^{0,\bullet}(X_t, E_t) \ar[d]^{\bar{\partial}_t} \\
  A^{0,\bullet+1}(X, E) \ar[r]^{\rho} & A^{0,\bullet+1}(X_t, E_t),   }
\]
Hence, we have
\begin{theorem}\label{thm-ext-iso}
The complex $(A^{0,\bullet}(X, E), \bar{\partial}_{\phi(t)} )$ and $(A^{0,\bullet}(X_t, E_t), \bar{\partial}_t)$ are isomorphic.
In particular, we have the isomorphism
\begin{equation}\label{ext-iso for coho}
H^{0,q}_{\bar{\partial}_{\phi(t)}}(X,E) \cong H^{0,q}_{\bar{\partial}_{t}}(X_t,E_t) ,\qquad\forall q \geq 0.
\end{equation}
\end{theorem}
\section{Deformations of holomorphic sections}
\label{deformation of holomorphic sections}

In this section, we consider a special case of the deformations of Dolbeault cohomology classes, i.e. the deformations of holomorphic sections. This will be the prototype of the more general cases that is studied in the following sections.

We will basically follow the terminology as given in \cite{Cat13}. Thus a \emph{deformation} of a compact complex space $X$ is a flat proper morphism $\varpi: (\mathcal{Y}, Y_{s_0}) \to (D, s_0)$ between connected complex spaces with $Y_{s_0}:=\varpi^{-1}(s_0)\cong X$ and a germ of deformation is called a \emph{small deformation}. It is also customary to call the fibers of $\varpi$ small deformations of $X$. We will also employ this convention frequently in this paper.

In this section and the following sections, we will make the following convention on our notations.  Let $\pi: (\mathcal{X}, X) \to (B, 0)$ be a small deformation of $X$ with each fiber of $\pi$ is a complex manifold, we will denote a holomorphic tensor bundle on $X$ by $E$ and for any $t\in B$, $E_t$ is the corresponding holomorphic tensor bundle on on $X_t=\pi^{-1}(t)$. We will always assume $X$ has been equipped with a fixed Hermitian metric.

Since the deformations of holomorphic sections can be constructed in a similar way as the Kuranishi family, we will first give a quick review of this notion.
\subsection{The Kuranishi family}
Let $X$ be a compact complex manifold equipped with an auxiliary Hermitian metric. The Kuranishi family of $X$ is the unique (up to isomorphism) semiuniversal small deformation $\pi: \mathcal{X}\to \mathcal{B}$ with $\pi^{-1}(0)=X_{0} \cong X$, where $0\in \mathcal{B}$. The Kuranishi family $\pi: (\mathcal{X}, X) \to (\mathcal{B}, 0)$ can be constructed as follows:

First, recall that in Hodge theory, we have the following operators defined on compact Hermitian manifolds:
\begin{itemize}
  \item $\Box:= \bar{\partial}\bar{\partial}^* + \bar{\partial}^*\bar{\partial}$, the $\bar{\partial}$-Laplacian operator;
  \item $\mathcal{H}: A^{0,q}(X, E)\to \mathcal{H}^{0,q}(X, E)$ be the projection operator to the $\bar{\partial}$-harmonic space $\mathcal{H}^{0,q}(X, E) = \ker \Box\cap A^{0,q}(X, E)$;
  \item $G:=$ the Green operator for $\Box$,
\end{itemize}
such that the following Hodge decomposition holds:
\begin{equation}\label{Hodge decomposition for par-bar}
Id = \mathcal{H} + \Box~G,\qquad\text{on}~A^{0,\bullet}(X, E).
\end{equation}
Let $\eta_1,\cdots, \eta_r$ be a basis of $\mathcal{H}^{0,1}(X, T^{1,0})$, where $\mathcal{H}^{0,1}(X, T^{1,0})\subset A^{0,1}(X, T^{1,0})$ is the harmonic space of $(0,1)$-forms with values in $T^{1,0}$, and $r= \dim \mathcal{H}^{0,1}(X, T^{1,0})= \dim H^1(X, T^{1,0})$. Let $\phi (t)=\sum_{\mu=1}^{\infty} \phi_{\mu}$, where
\[
\phi_{\mu}= \sum_{\nu_1+\cdots+\nu_r=\mu}^{\infty} \phi_{\nu_1 \cdots \nu_r}t_1^{\nu_1}\cdots t_r^{\nu_r}
\]
is the homogeneous term of degree $\mu$ and each $\phi_{\nu_1 \cdots \nu_r}\in A^{0,1}(X, T^{1,0})$. It can be shown that the following equation\footnote{Here and throughout this paper, $G$ will always denote the $\bar{\partial}$-Green operator unless otherwise stated.}
\begin{equation}
\left\{
\begin{array}{ll}
\phi (t)= \phi_1 + \frac{1}{2}\bar{\partial}^*G[\phi (t),\phi (t)], & \\
\phi_1 =\sum_{\nu=1}^{m}\eta_{\nu}t_{\nu}, & \\
\end{array} \right.
\end{equation}
has a unique power series solution $\phi (t)=\sum_{\mu=1}^{\infty} \phi_{\mu}$ which converges for $|t|$ small, and the Maurer-Cartan equation
\begin{equation}\label{Maurer-Cartan eq}
\bar{\partial}\phi (t)- \frac{1}{2}[\phi (t),\phi (t)]=0
\end{equation}
is satisfied if and only if $\mathcal{H}[\phi (t),\phi (t)]=0$, where $\mathcal{H}: A^{0,2}(X, T^{1,0})\to \mathcal{H}^{0,2}(X, T^{1,0})$ is the projection to harmonic space. For $\epsilon>0$ small, set
\[
\mathcal{B}:=\{t\in \mathbb{C}^r: |t|<\epsilon~\text{and}~\mathcal{H}[\phi (t),\phi (t)]=0 \}.
\]
Note that the Kuranishi space $\mathcal{B}$ is a complex analytic subset around the origin in $\mathbb{C}^r$ and $0\in \mathcal{B}$ since $\phi (0)=0$ . We can put a natural complex analytic space structure $\mathcal{X}$ on $X \times \mathcal{B}$ such that 1. $\pi: \mathcal{X}\to \mathcal{B}$ is a flat proper morphism between complex analytic spaces; 2. each fiber $X_t =\pi^{-1}(t)$ is the complex manifold obtained by endowing $X_t$
with the complex structure defined by $\phi (t)$ and $X_0=X$. For more information about the Kuranishi family, we refer the reader to \cite{Cat85,Cat13,Bal10,MK71}.

We will need the following:
\begin{lemma}\cite[pp.\,160]{MK71} \label{basic estimate}
 Let $\varphi\in A^{0,q}(X,E)$ and $G: A^{0,q}(X,E)\to A^{0,q}(X,E)$ be the Green operator, then for $k\geq 2$ we have
\begin{equation}
\|G\varphi\|_{k+\alpha}\leq C \|\varphi\|_{k-2+\alpha} ,
\end{equation}
where $C>0$ is independent of $\varphi$ and $\|\cdot\|_{k+\alpha}$ is the H\"older norm.
\end{lemma}

\begin{lemma}\label{lem-series-convergence}
Let $\sum_{k=0}^\infty x_kt^k$ and $\sum_{k=1}^\infty y_kt^k$ be two power series with nonnegative coefficients. Assume $\sum_{k=1}^\infty y_kt^k$ converges for small $t\in \mathbb{R}$ and for each $k>0$,
\[
x_k\leq C\sum_{j=1}^ky_jx_{k-j}
\]
holds for some constant $C>0$. Then $\sum_{k=0}^\infty x_kt^k$ converges for small $t\in\mathbb{R}$.
\end{lemma}
\begin{proof}By assumptions, we have $x_kt^k\leq C\sum_{j=1}^ky_jx_{k-j}t^k$ for each $k>0$ which implies for any $t>0$ and $m\geq 1$,
\[
\sum_{k=1}^m x_kt^k\leq C\sum_{k=1}^m\sum_{j=1}^ky_jx_{k-j}t^k\leq C\sum_{j=1}^my_jt^j\sum_{j=0}^{m-1}x_jt^j
\leq C\sum_{j=1}^my_jt^j\sum_{j=1}^mx_jt^j+x_0\sum_{j=1}^my_jt^j.
\]
Hence, for any small enough $t>0$ and $m\geq 1$,
\[
\sum_{k=1}^m x_kt^k\leq \frac{x_0\sum_{j=1}^my_jt^j}{1-C\sum_{j=1}^my_jt^j}<\infty.
\]
The conclusion then follows.
\end{proof}

\begin{theorem}\label{deformation of holomorphic sections: special case}
Let $\pi: (\mathcal{X}, X)\to (B,0)$ be a small deformation of a compact Hermitian manifold $(X,h)$ with Beltrami differentials $\phi(t)$. Assume $H^1(X,E)=0$ and $B$ is smooth.
Then given any $\sigma_0\in H^{0}(X,E)$, there exists a canonical family of holomorphic sections given by $\rho\sigma(t)\in H^{0}(X_t,E_t), \forall t\in B$ such that $\sigma(0)=\sigma_0$, $\sigma(t)=\sum_{k\geq 0} \sigma_k$ and
\begin{equation}
\forall k> 0,~~\sigma_{k}=\bar{\partial}^*G\sum_{i+j=k} \langle\phi_j | \sigma_{i} \rangle~~~~\in A^{0}(X,E).
\end{equation}
\end{theorem}

\begin{proof}The proof of this Theorem will be consist of three parts: existence, convergence and regularity. We may assume $B$ is a small polydisc with coordinates $t=(t_1,t_2,\cdots)$, then we may write $\phi=\phi (t)=\sum_{j=1}^{\infty} \phi_{j}$ as a power series in $t$ with each $\phi_{j}$ the homogenous term of degree $j$.

Existence: given $\sigma_0\in H^{0}(X,E)$, we want to find a power series $\sigma(t)=\sum_{k\geq 0} \sigma_k$ such that $\rho\sigma(t)\in A^{0}(X_t,E_t)$ is $\bar{\partial}_t$-closed which by Theorem \ref{ext formula for general holo tensor bundle} is equivalent to
\begin{equation}\label{Existence eq 1}
( \bar{\partial}-\langle\phi | ) \sigma(t)=0.
\end{equation}
By decomposing $\langle\phi | = \langle\phi (t)|$ and $\sigma(t)$ according to their degrees in $t$, \eqref{Existence eq 1} is reduced to the following system of equations :
\begin{equation}\label{Existence eq 2}
\left\{
\begin{array}{ll}
\bar{\partial}\sigma_k= \sum_{j=1}^{k}\langle\phi_j |\sigma_{k-j} \rangle, &~~k>0,  \\
\bar{\partial}\sigma_0=0. & \\
\end{array} \right.
\end{equation}
We will solve \eqref{Existence eq 2} by induction. Indeed, for $k=0$, \eqref{Existence eq 2} already has a solution $\sigma_0$. We therefore assume \eqref{Existence eq 2} can be solved for $k< N$ and our aim is to find $\sigma_{N}\in A^{0}(X, E)$ such that
\begin{equation}\label{Existence eq 3}
\bar{\partial}\sigma_{N}= \sum_{j=1}^{N}\langle\phi_j |\sigma_{N-j}\rangle.
\end{equation}
This equation can be solved by using Hodge theory after we have checked that
\begin{equation}\label{Existence eq 4}
\bar{\partial}\sum_{j=1}^{N}\langle\phi_j |\sigma_{N-j}\rangle=0.
\end{equation}
In fact, if \eqref{Existence eq 4} holds, then by Hodge theory, we have
\begin{equation}\label{Existence eq 5}
\sum_{j=1}^{N}\langle\phi_j |\sigma_{N-j}\rangle=\mathcal{H}\sum_{j=1}^{N}\langle\phi_j |\sigma_{N-j}\rangle+ \bar{\partial}\bar{\partial}^*G\sum_{j=1}^{N}\langle\phi_j |\sigma_{N-j}\rangle,
\end{equation}
where $\mathcal{H}: A^{0,1}(X, E)\to \mathcal{H}^1(X, E)$ is the projection to harmonic space and $G:A^{0,1}(X, E)\to A^{0,1}(X, E)$ is the Green operator. But our assumption $\mathcal{H}^1(X,E)=0$ implies $\mathcal{H}\sum_{j=1}^{N}\langle\phi_j | \sigma_{N-j}\rangle=0$. Hence we find a canonical solution of \eqref{Existence eq 3} given by
\begin{equation}
\sigma_{N}=\bar{\partial}^*G\sum_{j=1}^{N}\langle\phi_j |\sigma_{N-j}\rangle.
\end{equation}
Now we check \eqref{Existence eq 4}. Note that since $B$ is smooth, the Maurer-Cartan equation \eqref{Maurer-Cartan eq} is now reduced to
\begin{equation}\label{Maurer-Cartan eq decomposed}
\bar{\partial}\phi_k= \frac{1}{2}\sum_{j=1}^{k}[\phi_j,\phi_{k-j}],~~~~~\forall k>0.
\end{equation}
By \eqref{DGLA-mod-compatible},\eqref{phi-psi}, \eqref{Existence eq 2} and \eqref{Maurer-Cartan eq decomposed}, we have
\begin{align*}
&\bar{\partial}\sum_{j=1}^{N}\langle\phi_j |\sigma_{N-j}\rangle\\
=&\sum_{j=1}^{N} \left( -\langle\phi_j | \bar{\partial}\sigma_{N-j}\rangle + \langle\bar{\partial}\phi_j | \sigma_{N-j} \rangle \right)\\
=&\sum_{j=1}^{N} \left( -\langle\phi_j | \sum_{i=1}^{N-j}\langle\phi_i | \sigma_{N-j-i}\rangle +
\frac{1}{2}\langle \sum_{i=1}^{j}[\phi_i, \phi_{j-i}] | \sigma_{N-j}\rangle \right)\\
=&-\sum_{j=1}^{N}\sum_{i=1}^{N-j} \langle\phi_j | \langle\phi_i |\sigma_{N-j-i} \rangle \rangle +
\frac{1}{2} \left(\sum_{j=1}^{N}\sum_{i=1}^{j}  \langle\phi_i |  \langle\phi_{j-i} | \sigma_{N-j} \rangle \rangle
+\sum_{j=1}^{N}\sum_{i=1}^{j} \langle\phi_{j-i} | \langle\phi_{i} | \sigma_{N-j} \rangle \rangle \right)\\
=&-\sum_{i+j+k=N} \langle\phi_{i} | \langle\phi_j | \sigma_{k} \rangle \rangle +
\frac{1}{2} \left(\sum_{a+b+c=N} \langle\phi_a | \langle\phi_b  |\sigma_{c} \rangle \rangle
+\sum_{r+s+t=N} \langle\phi_{r} | \langle\phi_{s} | \sigma_{t} \rangle \rangle \right)\\
=&0.
\end{align*}
The proof of existence is complete.

Convergence: we note that for any $\psi\in A^{0,1}(X, \Theta)$, $\sigma \in A^{0}(X,E)$ and $\varphi \in A^{1}(X,E)$, it follows from the definition of the H\"older norm $\|\cdot\|_{k+\alpha}$ that
\begin{equation}
 \| \langle \psi | \sigma \rangle  \|_{k+\alpha}\leq C_1\|\psi\|_{k+1+\alpha} \|\sigma\|_{k+1+\alpha},
\end{equation}
and
\begin{equation}
 \|\bar{\partial}^*\varphi\|_{k+\alpha}\leq C_2 \|\varphi\|_{k+1+\alpha},
\end{equation}
where $C_1$ is independent of $\psi, \sigma$ and $C_2>0$ is independent of $\varphi$. It follows from Lemma \ref{basic estimate} that
\begin{equation}
\|\sigma_{j}\|_{k+\alpha}=\|\bar{\partial}^*G\sum_{a+b=j} \langle\phi_a | \sigma_{b}  \rangle \|_{k+\alpha} \leq C_2CC_1 \sum_{a+b=j} \|\phi_{a}\|_{k+\alpha} \|\sigma_{b}\|_{k+\alpha}.
\end{equation}
Without loss of generality, we may assume $\dim B=1$. In this case, write $\sigma(t)=\sum_{j\geq 0} \sigma_jt^j$ and $\phi(t)=\sum_{j\geq 1} \phi_jt^j$. Then the convergence of $\sigma(t)=\sum_{j\geq 0} \sigma_j$ follows from Lemma \ref{lem-series-convergence} by letting $x_j=\|\sigma_{j}\|_{k+\alpha}$ and $y_j=\|\phi_{j}\|_{k+\alpha}$.

Regularity: from the convergence of $\sigma(t)=\sum_{j\geq 0} \sigma_j$ in the H\"older norm $\|\cdot\|_{k+\alpha}$ for $k\geq 2$ we know that $\sigma(t)$ is at least $C^2$ on $X$. The smoothness of $\sigma(t)$ follows from the hypoellipticity of the operator $\bar{\partial}_t$. Indeed, $\rho\sigma(t)$ is a harmonic section since it is both $\bar{\partial}_t$-closed and $\bar{\partial}_t^*$-closed.
\end{proof}

Note that for any given $\sigma_0\in H^{0}(X,E)$, the power series $\sigma (t)=\sum_{k} \sigma_k$ given by $\sigma_{k}=\bar{\partial}^*G\sum_{i+j=k} \langle\phi_j | \sigma_{i} \rangle$ always converges even without the condition $H^1(X,E)=0$. We will see (in Theorem \ref{Deformation of Dolbeault cohomology classes: special case}) that for the conclusion of Theorem \ref{deformation of holomorphic sections: special case} to hold, the smoothness of $B$ is not necessary.

\section{Deformation of Dolbeault cohomology classes}
\label{Deformation of Dolbeault cohomology classes}
Let $\pi: (\mathcal{X}, X) \to (B, 0)$ be a small deformation of $X$ such that for each $t\in B$ the complex structure on $X_t$ is represented by Beltrami differential $\phi(t)$ and $E$ a holomorphic tensor bundle on $X$. We may always regard $B$ as an analytic subset in a small polydisc with coordinates $t=(t_1,t_2,\cdots)$ and $0\in B$. In this case, we may write $\phi(t)=\sum_{j=1}^{\infty} \phi_{j}$ as a power series in $t$ with each $\phi_{j}$ the homogenous term of degree $j$. Note that the Maurer-Cartan equation $\bar{\partial}\phi (t)- \frac{1}{2}[\phi (t),\phi (t)]=0$ holds only for those $t\in B$.
\subsection{The case $H^{q+1}(X,E)=0$}
We first show the following
\begin{proposition}\label{prop-unique-power-seris-solution}
$1$. $\forall \sigma \in A^{0,q}(X,E)$, if $\bar{\partial}_{\phi(t)} \sigma = \bar{\partial}\sigma - \langle\phi(t) | \sigma \rangle =0$ and $\bar{\partial}^*\sigma=0$, then we must have
\[
\sigma = \mathcal{H}\sigma + \bar{\partial}^*G \langle\phi(t) | \sigma \rangle ,
\]
where $\mathcal{H}: A^{0,q}(X, E)\to \mathcal{H}^{0,q}(X, E)$ is the projection operator to harmonic space.\\
$2$. For any fixed $\sigma_0\in \mathcal{H}^{0,q}(X,E)$, the equation
\begin{equation}\label{Kuranishi eq}
\sigma=\sigma_0 + \bar{\partial}^*G \langle\phi(t) | \sigma \rangle,
\end{equation}
has an unique solution given by $\sigma=\sigma(t)=\sum_{k} \sigma_k$ and $\sigma_{k}=\bar{\partial}^*G\sum_{i+j=k} \langle\phi_i | \sigma_{j} \rangle\in A^{0,q}(X,E) $ which converges for $|t|$ small.
\end{proposition}

\begin{proof}The first assertion follows from the Hodge decomposition:
\[
\sigma = \mathcal{H}\sigma + G\bar{\partial}^*\bar{\partial} \sigma + G\bar{\partial}\bar{\partial}^* \sigma
       = \mathcal{H}\sigma + G\bar{\partial}^* \langle\phi(t) | \sigma \rangle.
\]
For the second assertion, substitute $\sigma=\sigma(t)=\sum_{k} \sigma_k$ in \eqref{Kuranishi eq}, we have
\begin{equation}\label{formal solutions}
\left\{
\begin{array}{ll}
\sigma_1 &= \bar{\partial}^*G \langle\phi_1 | \sigma_0 \rangle,  \\
\sigma_2 &= \bar{\partial}^*G ( \langle\phi_2 | \sigma_0 \rangle + \langle\phi_1 | \sigma_1 \rangle ),\\
         &\cdots,  \\
\sigma_k &= \bar{\partial}^*G\sum_{i+j=k} \langle\phi_i | \sigma_j \rangle ,~~\forall k > 0.  \\
\end{array} \right.
\end{equation}
The convergence of $\sigma(t)$ follows from the same arguments as in the proof of Theorem \ref{deformation of holomorphic sections: special case}.

Uniqueness: Let $\sigma$ and $\sigma'$ be two solutions to $\sigma=\sigma_0 + \bar{\partial}^*G \langle\phi(t) | \sigma \rangle$ and set $\tau=\sigma-\sigma'$. Then $\tau=\bar{\partial}^*G \langle \phi (t) | \tau\rangle$ and so by Lemma \ref{basic estimate}, we have
\begin{equation}
\|\tau\|_{k+\alpha}\leq c \|\phi (t)\|_{k+\alpha}  \|\tau\|_{k+\alpha},
\end{equation}
for some constant $c>0$. When $|t|$ is sufficiently small, $\|\phi (t)\|_{k+\alpha}$ is also small. Hence we must have $\tau=0$.
\end{proof}

\begin{proposition}\label{prop-holo=harmonic part vanish}
Let $\sigma$ be a solution of the equation \eqref{Kuranishi eq} in Proposition \ref{prop-unique-power-seris-solution}.
Then $\forall t\in \mathcal{B}$, we have
\begin{equation}\label{eq-equivalence-harmonic part vanish}
\bar{\partial}\sigma = \langle\phi(t) | \sigma \rangle  \Leftrightarrow \mathcal{H} \langle\phi(t) | \sigma \rangle =0,
\end{equation}
where $\mathcal{H}: A^{0,q+1}(X, E)\to \mathcal{H}^{0,q+1}(X, E)$ is the projection operator to harmonic space. In particular, we have $\bar{\partial}\langle\phi(t) | \sigma \rangle=0$.
\end{proposition}

\begin{proof}First, if $\bar{\partial}\sigma = \langle\phi(t) | \sigma \rangle$, then it is clear that $\mathcal{H} \langle\phi(t) | \sigma \rangle =0$. Conversely, assume $\mathcal{H} \langle\phi(t) | \sigma \rangle =0$ and set
\[ \psi(t) := \bar{\partial}\sigma - \langle\phi(t) | \sigma \rangle.
\]
By our assumption, the Hodge decomposition of $\langle\phi(t) | \sigma \rangle$ now takes the form
\[
\langle\phi(t) | \sigma \rangle = \bar{\partial}\bar{\partial}^*G \langle\phi(t) | \sigma \rangle+\bar{\partial}^*\bar{\partial}G \langle\phi(t) | \sigma \rangle.
\]
It follows from \eqref{phi-phi} and \eqref{Kuranishi eq} that
\begin{align*}
\psi(t)
=& \bar{\partial}\bar{\partial}^*G \langle\phi(t) | \sigma \rangle - \langle\phi(t) | \sigma \rangle \\
=& -\bar{\partial}^*G\bar{\partial} \langle\phi(t) | \sigma \rangle \\
=& -\bar{\partial}^*G\left( \langle \bar{\partial}\phi(t) | \sigma \rangle - \langle\phi(t) | \bar{\partial} \sigma \rangle \right) \\
=& -\bar{\partial}^*G\left( \langle \bar{\partial}\phi(t) | \sigma \rangle - \langle\phi(t) | \psi(t) \rangle -  \langle\phi(t) |^2 \sigma \right)\\
=& -\bar{\partial}^*G\left( \langle \bar{\partial}\phi(t) - \frac{1}{2}[\phi(t) ,\phi(t) ] | \sigma \rangle - \langle \phi(t) | \psi(t) \rangle \right)\\
=& \bar{\partial}^*G \langle\phi(t) | \psi(t) \rangle .
\end{align*}
It follows that
\begin{equation}
\|\psi(t)\|_{k+\alpha}=\|\bar{\partial}^*G \langle\phi(t) | \psi(t) \rangle\|_{k+\alpha} \leq C \|\phi(t)\|_{k+\alpha}\|\psi(t)\|_{k+\alpha},
\end{equation}
where $C$ depends only on $k$ and $\alpha$. Now choose $|t|$ so small that $C \|\phi(t)\|_{k+\alpha}<1$, then we get $\|\psi(t)\|_{k+\alpha}< \|\psi(t)\|_{k+\alpha}$ which is a contradiction if $\psi(t)\neq0$. So we must have $\psi(t)=0$ whenever $|t|$ is small enough.

At last, since \eqref{Kuranishi eq} implies $\bar{\partial}\sigma=\bar{\partial}\bar{\partial}^*G\langle\phi(t) | \sigma \rangle$, we have
\begin{align*}
\langle\phi(t) | \sigma \rangle=&
\mathcal{H} \langle\phi(t) | \sigma \rangle+\bar{\partial}\bar{\partial}^*G\langle\phi(t) | \sigma \rangle+\bar{\partial}^*\bar{\partial}G\langle\phi(t) | \sigma \rangle\\
=& \mathcal{H} \langle\phi(t) | \sigma \rangle+\bar{\partial}\sigma +\bar{\partial}^*\bar{\partial}G\langle\phi(t) | \sigma \rangle.
\end{align*}
We see that \eqref{eq-equivalence-harmonic part vanish} holds if and only if $\bar{\partial}^*\bar{\partial}G\langle\phi(t) | \sigma \rangle=0$ which is equivalent to $\bar{\partial}\langle\phi(t) | \sigma \rangle=0$.
\end{proof}

\begin{theorem}\label{Deformation of Dolbeault cohomology classes: special case}
Let $\pi: (\mathcal{X}, X)\to (B,0)$ be a small deformation of a compact Hermitian manifold $(X,h)$ with Beltrami differentials $\phi(t)$. Assume that $H^{q+1}(X,E)=0$.
Then given any $\sigma_0\in \mathcal{H}^{0,q}(X,E)$, there exists a canonical family of $\bar{\partial}_{\phi}$-closed $E$-valued $(0,q)$-forms given by $\sigma(t)\in A^{0,q}(X,E), \forall t\in B$ such that $\sigma(0)=\sigma_0$, $\sigma(t)=\sum_{k} \sigma_k$ and
\begin{equation}
\forall k> 0,~~\sigma_{k}=\bar{\partial}^*G\sum_{i+j=k} \langle\phi_i | \sigma_j \rangle~~~~\in A^{0,q}(X,E).
\end{equation}
\end{theorem}

\begin{proof}It follows from the assumption $H^{q+1}(X,E)=0$ and Proposition \ref{prop-holo=harmonic part vanish} that the power series $\sigma(t)=\sum_{k} \sigma_k$ is $\bar{\partial}_{\phi}$-closed. We are left to show that $\sigma(t)$ is smooth. Indeed, from the convergence of $\sigma(t)=\sum_{k} \sigma_k$ in the H\"older norm $\|\cdot\|_{k+\alpha}$ for $k\geq 2$ we know that $\sigma(t)$ is at least $C^2$ on $X$. Let $\Box$ be the $\bar{\partial}$-Laplacian operator, it follows from \eqref{Kuranishi eq} and Hodge theory that
\begin{align*}
\Box \sigma(t)
=& \Box\bar{\partial}^*G \langle\phi | \sigma (t)\rangle\\
=& \bar{\partial}^*\Box G \langle\phi | \sigma (t)\rangle\\
=& \bar{\partial}^* \langle\phi | \sigma (t)\rangle - \bar{\partial}^*\mathcal{H} \langle\phi | \sigma (t)\rangle\\
=& \bar{\partial}^* \langle\phi | \sigma (t)\rangle .
\end{align*}
We see that $\sigma(t)$ satisfies the following equation:
\begin{equation}\label{regularity eq}
\Box \sigma(t) - \bar{\partial}^* \langle\phi(t) | \sigma (t)\rangle = 0,
\end{equation}
which is a (small perturbation of a) standard elliptic equation for small $t$. Hence we may conclude that $\sigma(t)$ is smooth on $X$ for any fixed $t\in B$.
\end{proof}
\begin{remark}
In Theorem \ref{Deformation of Dolbeault cohomology classes: special case}, we do not need to assume $B$ is smooth. Because when $B$ is not smooth, we may assume $B$ is an analytic subset in a small polydisc with $0\in B$. The above proof shows we can construct a convergent power series $\sigma(t)$ with constant term equal to $\sigma_0$ and such that for any fixed $t\in B$, $\sigma(t)$ is a $\bar{\partial}_{\phi}$-closed (by Proposition \ref{prop-holo=harmonic part vanish}) smooth $E$-valued $(0,q)$-form on $X$. Note that the convergent power series $\sigma(t)$ (whose explicit form is suggested by the proof of Theorem \ref{deformation of holomorphic sections: special case} and Proposition \ref{prop-unique-power-seris-solution}) can always be defined with or without the smoothness of $B$. But in Theorem \ref{deformation of holomorphic sections: special case} the smoothness of $B$ is assumed since in the proof we have constructed the convergent power series $\sigma(t)=\sum_{j\geq 0} \sigma_j$ inductively. For this purpose, we must use equation \eqref{Maurer-Cartan eq decomposed} which is a consequence of the smoothness of $B$.
\end{remark}

\subsection{Deformation of Dolbeault cohomology classes: the general case}
First, we observe that the operator $\bar{\partial}_{\phi} = \bar{\partial} - \langle\phi|$ is, like the operator $\bar{\partial}$, an elliptic differential operator of degree $1$. This follows from the fact that $\bar{\partial}$ is elliptic and $\phi=\phi(t)$ represents a small deformation. As a result, if we set
\begin{itemize}
  \item $\Box_{\phi}:= \bar{\partial}_{\phi}\bar{\partial}_{\phi}^* + \bar{\partial}_{\phi}^*\bar{\partial}_{\phi}$, the $\bar{\partial}_{\phi}$-Laplacian operator;
  \item $\mathcal{H}_{\phi}: A^{0,q}(X, E)\to \mathcal{H}_{\phi}^{0,q}(X, E)$ be the projection operator to the $\bar{\partial}_{\phi}$-harmonic space $\mathcal{H}_{\phi}^{0,q}(X, E) = \ker \Box_{\phi}\cap A^{0,q}(X, E)$;
  \item $G_{\phi}:=$ the Green operator for $\Box_{\phi}$.
\end{itemize}
then the following Hodge decomposition holds~\cite[Chap.\,VI]{Dem12}:
\begin{equation}\label{Hodge decomposition for par-bar-phi}
I = \mathcal{H}_{\phi} + \Box_{\phi}G_{\phi}.
\end{equation}
Notice that while $\bar{\partial}_{\phi(t)}$ depends holomorphically on $t$, its adjoint operator $\bar{\partial}_{\phi(t)}^*$ depends anti-holomorphically on $t$. So we have $\Box_{\phi(t)}$ depends smoothly on $t$. But this does not hold in general for $G_{\phi(t)}$ or $\mathcal{H}_{\phi(t)}$. In fact, the continuity in $t$ of $G_{\phi(t)}$ and $\mathcal{H}_{\phi(t)}$, as operators on $A^{0,q}(X, E)$, is equivalent to the condition that $\dim \mathcal{H}_{\phi(t)}^{0,q}(X, E)$ is independent of $t$, see Theorem $4.5$ of \cite[pp.\,178]{MK71}.
It can be shown as usual that
\[
\bar{\partial}_{\phi}\mathcal{H}_{\phi}=\mathcal{H}_{\phi}\bar{\partial}_{\phi}=\bar{\partial}_{\phi}^*\mathcal{H}_{\phi}=\mathcal{H}_{\phi}\bar{\partial}_{\phi}^*
=G_{\phi}\mathcal{H}_{\phi}=\mathcal{H}_{\phi}G_{\phi}=[\bar{\partial}_{\phi},G_{\phi}]=[\bar{\partial}_{\phi}^*,G_{\phi}]=0,
\]
and
\[
\| \langle \phi | \sigma \rangle \|_{k+\alpha} \leq c_{k+\alpha} \| \phi \|_{k+\alpha+1} \| \sigma \|_{k+\alpha+1},
\]
for any $\phi\in A^{0,l}(X, T_{X}^{1,0})$ and $\sigma\in A^{0,q}(X, E)$.

\begin{proposition}\label{partial-bar*-closed-rep of deformed-Dol-coho}
The following natural homomorphism induced by inclusion $\ker\bar{\partial}^*\cap \ker\bar{\partial}_{\phi(t)}\subseteq\ker\bar{\partial}_{\phi(t)}$ is an isomorphism\footnote{To save space, we have omitted $\cap A^{0,q}(X,E)$ which will be clear from the context.}:
\[
\frac{\ker\bar{\partial}^*\cap \ker\bar{\partial}_{\phi(t)}}{\ker\bar{\partial}^*\cap \Image\bar{\partial}_{\phi(t)}}\longrightarrow H_{\phi(t)}^{0,q}(X, E).
\]
\end{proposition}

\begin{proof}By Hodge theory for the deformed Laplacian we have
\[
(\ker\bar{\partial}_{\phi(t)})^\perp=\Image\bar{\partial}_{\phi(t)}^*
\quad\text{and}\quad
(\ker\bar{\partial}^*)^\perp=\Image\bar{\partial}
\]
which implies the following orthogonal direct sum decomposition
\[
A^{0,q}(X,E)=\left(\ker\bar{\partial}_{\phi(t)}\cap \ker\bar{\partial}^*\right)\oplus (\Image\bar{\partial}_{\phi(t)}^*+\Image\bar{\partial}),
\]
from which we deduce that
\[
\ker\bar{\partial}_{\phi(t)}=\left(\ker\bar{\partial}_{\phi(t)}\cap \ker\bar{\partial}^*\right)\oplus
\left(\ker\bar{\partial}_{\phi(t)}\cap (\Image\bar{\partial}+\Image\bar{\partial}_{\phi(t)}^*)\right).
\]
Similarly, we have
\begin{align*}
\Image\bar{\partial}_{\phi(t)}&=
\left(\Image\bar{\partial}_{\phi(t)}\cap \ker\bar{\partial}^*\right)\oplus
\left(\Image\bar{\partial}_{\phi(t)}\cap (\Image\bar{\partial}_{\phi(t)}\cap \ker\bar{\partial}^*)^\perp\right)
\\
&= \left(\Image\bar{\partial}_{\phi(t)}\cap \ker\bar{\partial}^*\right)\oplus
\left(\Image\bar{\partial}_{\phi(t)}\cap (\ker\bar{\partial}_{\phi(t)}^*+\Image\bar{\partial})\right),
\end{align*}
and
\begin{align}
\ker\bar{\partial}_{\phi(t)}&\cap (\Image\bar{\partial}+\Image\bar{\partial}_{\phi(t)}^*)=
\left(\Image\bar{\partial}_{\phi(t)}\cap(\Image\bar{\partial}+\Image\bar{\partial}_{\phi(t)}^*)\right)\oplus \label{align-1}\\
&\ker\bar{\partial}_{\phi(t)}\cap (\Image\bar{\partial}+\Image\bar{\partial}_{\phi(t)}^*)\cap(\ker\bar{\partial}_{\phi(t)}^*+\ker\bar{\partial}^*\cap\ker\bar{\partial}_{\phi(t)} )\notag
\\
\Image\bar{\partial}_{\phi(t)}&\cap (\ker\bar{\partial}_{\phi(t)}^*+\Image\bar{\partial})=
\left(\Image\bar{\partial}_{\phi(t)}\cap(\Image\bar{\partial}+\Image\bar{\partial}_{\phi(t)}^*)\right)\oplus\label{align-2}\\
&\Image\bar{\partial}_{\phi(t)}\cap (\ker\bar{\partial}_{\phi(t)}^*+\Image\bar{\partial})\cap(\ker\bar{\partial}_{\phi(t)}^*+\ker\bar{\partial}^*\cap\ker\bar{\partial}_{\phi(t)} ).\notag
\end{align}
It follows from \eqref{align-1} that
\begin{align*}
&\frac{\ker\bar{\partial}_{\phi(t)}\cap\left(\Image\bar{\partial}+
\Image\bar{\partial}_{\phi(t)}^*\right)}{\Image\bar{\partial}_{\phi(t)}\cap\left(\Image\bar{\partial}+\Image\bar{\partial}_{\phi(t)}^*\right)}\\
\cong& \ker\bar{\partial}_{\phi(t)}\cap (\Image\bar{\partial}+\Image\bar{\partial}_{\phi(t)}^*)\cap(\ker\bar{\partial}_{\phi(t)}^*+\ker\bar{\partial}_{\phi(t)}\cap \ker\bar{\partial}^*)\\
=& (\Image\bar{\partial}+\Image\bar{\partial}_{\phi(t)}^*)\cap(\ker\bar{\partial}_{\phi(t)}\cap\ker\bar{\partial}_{\phi(t)}^*+\ker\bar{\partial}_{\phi(t)}\cap \ker\bar{\partial}^*),
\end{align*}
and from \eqref{align-2} that
\begin{align*}
&\frac{\Image\bar{\partial}_{\phi(t)}\cap\left(\Image\bar{\partial}
+\ker\bar{\partial}_{\phi(t)}^*\right)}{\Image\bar{\partial}_{\phi(t)}\cap\left(\Image\bar{\partial}+\Image\bar{\partial}_{\phi(t)}^*\right)}\\
\cong& \Image\bar{\partial}_{\phi(t)}\cap (\Image\bar{\partial}+\ker\bar{\partial}_{\phi(t)}^*)\cap(\ker\bar{\partial}_{\phi(t)}^*+\ker\bar{\partial}_{\phi(t)}\cap \ker\bar{\partial}^*)\\
=& \Image\bar{\partial}_{\phi(t)}\cap (\Image\bar{\partial}+\ker\bar{\partial}_{\phi(t)}^*)\cap(\ker\bar{\partial}_{\phi(t)}\cap\ker\bar{\partial}_{\phi(t)}^*+\ker\bar{\partial}_{\phi(t)}\cap \ker\bar{\partial}^*),
\end{align*}
which are finite dimensional. Indeed, both $\ker\bar{\partial}_{\phi(t)}\cap\ker\bar{\partial}_{\phi(t)}^*=\ker\Box_{\phi(t)}$ and $\ker\bar{\partial}_{\phi(t)}\cap \ker\bar{\partial}^*=\ker(\bar{\partial}_{\phi(t)}^*\bar{\partial}_{\phi(t)}+\bar{\partial}\bar{\partial}^*)$ are the null space of elliptic operators.
It follows that
\begin{align*}
&\dim H_{\phi(t)}^{0,q}(X, E)+\dim\frac{\Image\bar{\partial}_{\phi(t)}\cap\left(\Image\bar{\partial}
+\ker\bar{\partial}_{\phi(t)}^*\right)}{\Image\bar{\partial}_{\phi(t)}\cap\left(\Image\bar{\partial}+\Image\bar{\partial}_{\phi(t)}^*\right)}\\
=&\dim \frac{\ker\bar{\partial}_{\phi(t)}\cap\left(\Image\bar{\partial}+
\Image\bar{\partial}_{\phi(t)}^*\right)}{\Image\bar{\partial}_{\phi(t)}\cap\left(\Image\bar{\partial}+\Image\bar{\partial}_{\phi(t)}^*\right)}
+\dim\frac{\ker\bar{\partial}_{\phi(t)}\cap \ker\bar{\partial}^*}{\Image\bar{\partial}_{\phi(t)}\cap \ker\bar{\partial}^*}.
\end{align*}
For any $x\in \Image\bar{\partial}$, there exists unique $y\in\ker\bar{\partial}_{\phi(t)}$ and unique $z\in\Image\bar{\partial}_{\phi(t)}^*$ such that $x=y+z$. This defines a surjective homomorphism
\[
\Image\bar{\partial}\longrightarrow\ker\bar{\partial}_{\phi(t)}\cap\left(\Image\bar{\partial}+
\Image\bar{\partial}_{\phi(t)}^*\right):x\longmapsto y,
\]
with $\Image\bar{\partial}\cap\Image\bar{\partial}_{\phi(t)}^*$ as kernel. Repeating this argument, we get
\[
\Image\bar{\partial}_{\phi(t)}\cap\left(\Image\bar{\partial}+\Image\bar{\partial}_{\phi(t)}^*\right)\cong \frac{\Image\bar{\partial}\cap\Image\square_{\phi(t)}}{\Image\bar{\partial}\cap\Image\bar{\partial}_{\phi(t)}^*},
\]
and
\[
\Image\bar{\partial}_{\phi(t)}\cap\left(\Image\bar{\partial}
+\ker\bar{\partial}_{\phi(t)}^*\right)\cong \frac{\Image\bar{\partial}}{\Image\bar{\partial}\cap\ker\bar{\partial}_{\phi(t)}^*}.
\]
Hence we have
\begin{align*}
&\dim \frac{\ker\bar{\partial}_{\phi(t)}\cap\left(\Image\bar{\partial}+
\Image\bar{\partial}_{\phi(t)}^*\right)}{\Image\bar{\partial}_{\phi(t)}\cap\left(\Image\bar{\partial}+\Image\bar{\partial}_{\phi(t)}^*\right)}
-\dim\frac{\Image\bar{\partial}_{\phi(t)}\cap\left(\Image\bar{\partial}
+\ker\bar{\partial}_{\phi(t)}^*\right)}{\Image\bar{\partial}_{\phi(t)}\cap\left(\Image\bar{\partial}+\Image\bar{\partial}_{\phi(t)}^*\right)}\\
=&\dim\frac{\Image\bar{\partial}}{\Image\bar{\partial}\cap\Image\bar{\partial}_{\phi(t)}^*}/
\frac{\Image\bar{\partial}\cap\Image\square_{\phi(t)}}{\Image\bar{\partial}\cap\Image\bar{\partial}_{\phi(t)}^*}
-\dim \frac{\Image\bar{\partial}}{\Image\bar{\partial}\cap\ker\bar{\partial}_{\phi(t)}^*}/
\frac{\Image\bar{\partial}\cap\Image\square_{\phi(t)}}{\Image\bar{\partial}\cap\Image\bar{\partial}_{\phi(t)}^*}\\
=&\dim \frac{\Image\bar{\partial}\cap\ker\bar{\partial}_{\phi(t)}^*}{\Image\bar{\partial}\cap\Image\bar{\partial}_{\phi(t)}^*}.
\end{align*}
It follows that
\[
\dim H_{\phi(t)}^{0,q}(X, E)= \dim\frac{\Image\bar{\partial}\cap\ker\bar{\partial}_{\phi(t)}^*}{\Image\bar{\partial}\cap\Image\bar{\partial}_{\phi(t)}^*}+\dim
\frac{\ker\bar{\partial}^*\cap \ker\bar{\partial}_{\phi(t)}}{\ker\bar{\partial}^*\cap \Image\bar{\partial}_{\phi(t)}}~.
\]
It is left to show $\Image\bar{\partial}\cap\ker\bar{\partial}_{\phi(t)}^*=0$. Let $\sigma\in\ker\bar{\partial}\cap\ker\bar{\partial}_{\phi(t)}^*$, by the same argument as in Proposition \ref{prop-unique-power-seris-solution} we have that $\sigma$ must be the unique solution of the equation
\[
\sigma=\sigma_0+\bar{\partial}G\langle\phi(t)|^*\sigma,\quad \sigma_0=\mathcal{H}\sigma,
\]
where $\langle\phi(t)|^*=\bar{\partial}^*-\bar{\partial}_{\phi(t)}^*$. Now if $\sigma\in\Image\bar{\partial}\cap\ker\bar{\partial}_{\phi(t)}^*$ then $\sigma$ satisfies $\sigma=\sigma_0+\bar{\partial}G\langle\phi(t)|^*\sigma$ with $\sigma_0=\mathcal{H}\sigma=0\Rightarrow\sigma=0$.
\end{proof}

\begin{proposition}\label{prop-kerdbar*-imagedbarphi}
$1$. For any fixed $t\in B$, the following homomorphism
\[
g_t:\ker\bar{\partial}\cap A^{0,q}(X,E)\longrightarrow \ker\bar{\partial}^*\cap\Image\bar{\partial}_{\phi(t)}\cap A^{0,q+1}(X,E):x_0\longmapsto \bar{\partial}_{\phi(t)}x(t),
\]
is surjective with $\ker g_t=\ker\bar{\partial}\cap\left(\ker\bar{\partial}_{\phi(t)}+\Image\bar{\partial}^*\right)\cap A^{0,q}(X,E)$, where $x(t)$ is the unique solution of $x(t)=x_0 + \bar{\partial}^*G \langle\phi(t) | x(t)\rangle$. Furthermore, we have
\[
\bar{\partial}_{\phi(t)}x(t)=0\Leftrightarrow \mathcal{H}\langle\phi(t) | x(t) \rangle =0,
\]
and $\bar{\partial}\langle\phi(t) | x(t) \rangle =0$.\\
$2$. Let $\hat{g}_t:\mathcal{H}^{0,q}(X,E)\longrightarrow \ker\bar{\partial}^*\cap\Image\bar{\partial}_{\phi(t)}\cap A^{0,q+1}(X,E)$ be the restriction of $g_t$ on $\mathcal{H}^{0,q}(X,E)$, then $\hat{g}_t$ is surjective with
$\ker \hat{g}_t=\mathcal{H}^{0,q}(X,E)\cap\left(\ker\bar{\partial}_{\phi(t)}+\Image\bar{\partial}^*\right)$. Moreover, we have
\[
\dim\mathcal{H}^{0,q}(X,E)=\dim\left(\ker\bar{\partial}^*\cap\ker\bar{\partial}_{\phi(t)}\right)^q+
\dim \left(\ker\bar{\partial}^*\cap\Image\bar{\partial}_{\phi(t)}\right)^{q+1},
\]
where $\left(\ker\bar{\partial}^*\cap\ker\bar{\partial}_{\phi(t)}\right)^q:=\ker\bar{\partial}^*\cap\ker\bar{\partial}_{\phi(t)}\cap A^{0,q}(X,E)$ and $\left(\ker\bar{\partial}^*\cap\Image\bar{\partial}_{\phi(t)}\right)^{q+1}:=\ker\bar{\partial}^*\cap\Image\bar{\partial}_{\phi(t)}\cap A^{0,q+1}(X,E)$.
\end{proposition}
\begin{proof}$1$. Let $x\in A^{0,q}(X,E)$, then by Hodge decomposition we have
\[
\bar{\partial}_{\phi(t)}x=\bar{\partial}x-\langle\phi(t) | x\rangle
=\bar{\partial}x-\mathcal{H}\langle\phi(t) | x\rangle-\bar{\partial}^*\bar{\partial}G\langle\phi(t) | x\rangle- \bar{\partial}\bar{\partial}^*G\langle\phi(t) | x\rangle,
\]
thus
\[
\bar{\partial}_{\phi(t)}x\in\ker\bar{\partial}^*\Leftrightarrow \bar{\partial}x-\bar{\partial}\bar{\partial}^*G\langle\phi(t) | x\rangle=0.
\]
Set $x_0=x-\bar{\partial}^*G\langle\phi(t) | x\rangle\in\ker\bar{\partial}$, then $x$ is a solution to the equation $x=x_0 + \bar{\partial}^*G \langle\phi(t) | x\rangle$ which is uniquely determined by $x_0$ in view of the proof of Proposition \ref{prop-unique-power-seris-solution}. In fact, we have $x=x(t)=\sum_{k} x_k$ with $x_{k}=\bar{\partial}^*G\sum_{i+j=k} \langle\phi_i | x_{j} \rangle,~\forall k> 0$. The last assertion follows from the same proof as in Proposition \ref{prop-holo=harmonic part vanish}.

It is left to show $\ker g_t=\ker\bar{\partial}\cap\left(\ker\bar{\partial}_{\phi(t)}+\Image\bar{\partial}^*\right)\cap A^{0,q}(X,E)$. In fact, it is clear that $\ker g_t\subseteq\ker\bar{\partial}\cap\left(\ker\bar{\partial}_{\phi(t)}+\Image\bar{\partial}^*\right)$. Conversely, let us consider the following surjective homomorphism
\[
\ker\bar{\partial}_{\phi(t)}\longrightarrow \ker\bar{\partial}\cap\left(\ker\bar{\partial}_{\phi(t)}+\Image\bar{\partial}^*\right):
x\longmapsto x-\bar{\partial}^*\bar{\partial}Gx=x-\bar{\partial}^*G\langle\phi(t) | x\rangle,
\]
whose kernel is $\ker\bar{\partial}_{\phi(t)}\cap\Image\bar{\partial}^*=0$ by Proposition \ref{prop-unique-power-seris-solution}. Its inverse is given by
\[
\ker\bar{\partial}\cap\left(\ker\bar{\partial}_{\phi(t)}+\Image\bar{\partial}^*\right)\longrightarrow \ker\bar{\partial}_{\phi(t)}:
x_0\longmapsto x(t),
\]
where $x(t)$ is the unique solution of $x(t)=x_0 + \bar{\partial}^*G \langle\phi(t) | x(t)\rangle$.
As a result, let $x_0\in\ker\bar{\partial}\cap\left(\ker\bar{\partial}_{\phi(t)}+\Image\bar{\partial}^*\right)$ and $x(t)$ be the unique solution of $x(t)=x_0 + \bar{\partial}^*G \langle\phi(t) | x(t)\rangle$, we must have $x(t)\in \ker\bar{\partial}_{\phi(t)}\Rightarrow x_0\in\ker g_t$. \\

$2$. It can be proved in the same way that $\ker \hat{g}_t=\mathcal{H}^{0,q}(X,E)\cap\left(\ker\bar{\partial}_{\phi(t)}+\Image\bar{\partial}^*\right)$. To show $\hat{g}_t$ is surjective it is enough to show
\[
\frac{\mathcal{H}^{0,q}(X,E)}{\mathcal{H}^{0,q}(X,E)\cap\left(\ker\bar{\partial}_{\phi(t)}+\Image\bar{\partial}^*\right)}\cong
\frac{\ker\bar{\partial}}{\ker\bar{\partial}\cap\left(\ker\bar{\partial}_{\phi(t)}+\Image\bar{\partial}^*\right)}.
\]
Indeed, we have
\begin{align*}
&\ker\bar{\partial}\cap\left(\ker\bar{\partial}_{\phi(t)}+\Image\bar{\partial}^*\right)\cong\ker\bar{\partial}_{\phi(t)}\\
=&\left(\ker\bar{\partial}_{\phi(t)}\cap (\Image\bar{\partial}+\Image\bar{\partial}_{\phi(t)}^*)\right)\oplus
\left(\ker\bar{\partial}_{\phi(t)}\cap \ker\bar{\partial}^*\right).
\end{align*}
Similarly, we have
\[
\left(\ker\bar{\partial}_{\phi(t)}\cap (\Image\bar{\partial}+\Image\bar{\partial}_{\phi(t)}^*)\right)\cong \Image\bar{\partial},
\]
and
\[
\ker\bar{\partial}_{\phi(t)}\cap \ker\bar{\partial}^*\cong\mathcal{H}^{0,q}(X,E)\cap\left(\ker\bar{\partial}_{\phi(t)}+\Image\bar{\partial}^*\right).
\]
The conclusion then follows since $\ker\bar{\partial}=\mathcal{H}^{0,q}(X,E)\oplus\Image\bar{\partial}$.
\end{proof}
\begin{remark}\label{rk-kerdbar*-imagedbarphi}
It follows from this Proposition that $\ker\bar{\partial}^*\cap\Image\bar{\partial}_{\phi(t)}\cap A^{0,\bullet}(X,E)=\Image\bar{\partial}_{\phi(t)}\cap \mathcal{H}^{0,\bullet}(X,E)$ because
$\bar{\partial}\bar{\partial}_{\phi(t)}x(t)=\bar{\partial}\Big(\bar{\partial}x(t)-\langle\phi(t) | x(t) \rangle \Big)=0$.
\end{remark}
\begin{definition}\label{def-V_t}
For any $t\in B$ and a vector subspace $V=\mathbb{C}\{ \sigma_0^1, \cdots, \sigma_0^N \}\subseteq \mathcal{H}^{0,q}(X,E)$, we set
\begin{align*}
V_{t}:=
&\Big\{ \sum_{l=1}^N a_l\sigma_0^l\in V \mid  (a_1, \cdots, a_N)\in \mathbb{C}^{N}~\text{s.t.}~\bar{\partial}_{\phi(t)}\sigma(t)=0~\text{where}~\sigma(t)~\\
&\text{is the unique solution of the equation}~\sigma (t)=\sum_{l=1}^N a_l\sigma_0^l+\bar{\partial}^*G\langle\phi(t) | \sigma(t) \rangle\Big\}.
\end{align*}
Note that $V_{t}$ consists of those vectors of the form $\sum_l a_l\sigma_0^l$ such that the coefficients $a_l$ satisfy the following linear equation:
\[
\sum_{l=1}^N a_l\bar{\partial}_{\phi(t)}\sigma^l(t) =0,\quad\text{or}\quad \sum_{l=1}^N a_l \mathcal{H}\langle\phi(t) | \sigma^l(t) \rangle =0,
\]
where $\sigma^l(t)= \sum_{k} \sigma_k^l$ with $\sigma_{k}^l=\bar{\partial}^*G\sum_{i+j=k} \langle\phi_{i} | \sigma_{j}^l \rangle,~\forall k\neq 0$. We see that $V_{t}$ is a vector subspace of $V$ and varies holomorphically with $t$. In particular, the set $\{t\in B\mid\dim V_{t}\geq k\}$ is an analytic subset of $B$ for any $k\in \mathbb{N}$.
\end{definition}

\begin{definition}\label{def-f_t-g_t}
In view of Proposition \ref{partial-bar*-closed-rep of deformed-Dol-coho}, we set
\begin{align*}
f_t: &V_{t} \longrightarrow \frac{\ker\bar{\partial}^*\cap\ker\bar{\partial}_{\phi(t)}}{\ker\bar{\partial}^*\cap \Image\bar{\partial}_{\phi(t)}}\cong H^{0,q}_{\bar{\partial}_{\phi(t)}}(X,E),\\
&\sigma_0\longmapsto \sigma(t)=\sum_{k} \sigma_k,~\text{where}~\sigma_{k}=\bar{\partial}^*G\sum_{i+j=k} \langle\phi_i | \sigma_{j} \rangle,~\forall k> 0.
\end{align*}
\end{definition}

\begin{proposition}\label{prop V-E-t}
If $V= \mathcal{H}^{0,q}(X,E)$, then $f_t$ is surjective.
\end{proposition}

\begin{proof}By Proposition \ref{prop-unique-power-seris-solution}, the map
\begin{align*}
\tilde{f}_t: &V_{t} \longrightarrow \ker\bar{\partial}^*\cap\ker\bar{\partial}_{\phi(t)}\cap A^{0,q}(X,E),\\
&\sigma_0\longmapsto \sigma(t)=\sum_{k} \sigma_k,~\text{where}~\sigma_{k}=\bar{\partial}^*G\sum_{i+j=k} \langle\phi_i | \sigma_{j} \rangle,~\forall k\neq 0,
\end{align*}
is an isomorphism. Indeed, the inverse map of $\tilde{f}_t$ is given by the harmonic projection $\mathcal{H}$.
\end{proof}
For $V= \mathcal{H}^{0,q}(X,E)$, we see that an immediate consequence of Proposition \ref{prop V-E-t} is
\[
H^{q}(X_t,E_t)\cong V_{t}/\ker f_t ,
\]
for any $t\in B$.

\begin{theorem}\label{Deformation of Dolbeault cohomology classes: general case}
Let $\pi: (\mathcal{X}, X)\to (B,0)$ be a small deformation of a compact Hermitian manifold $(X,h)$ with Beltrami differentials $\phi(t)$. Let $V=\mathbb{C}\{ \sigma_0^1, \cdots, \sigma_0^N \}$ be a linear subspace of $\mathcal{H}^{0,q}(X,E)$ and $\sigma^l(t)=\tilde{f}_t\sigma_{0}^l,~l=1, \cdots, N$. Define an analytic subset $B(V)$ of $B$ by
\[
B(V):=\{t\in B\mid \mathcal{H} \langle\phi(t) | \sigma^l(t) \rangle = 0, l=1, \cdots, N\},
\]
Then we have
\begin{equation}\label{B(V)-discription}
B(V)=\{t\in B\mid \dim V = \dim \Image f_t +\dim \ker f_t\}.
\end{equation}
In particular, we have
\begin{equation}\label{B'-discription}
B':=B(\mathcal{H}^{0,q}(X,E))=\{t\in B\mid \dim H^{q}(X,E)=\dim H^{q}(X_t,E_t)+\dim \ker f_t\} .
\end{equation}
\end{theorem}

\begin{proof}
First, let $\alpha_1, \cdots, \alpha_d$ be a basis of the harmonic space $\mathcal{H}^{q+1}(X,E)$, then $\forall l=1, \cdots, N$, we have
\[
\mathcal{H} \langle\phi(t) | \sigma^l(t) \rangle = 0 \Leftrightarrow < \langle\phi(t) | \sigma^l(t) \rangle ,\alpha_q>=0,~\forall q=1, \cdots, d,
\]
where $<\cdot,\cdot>$ is the inner product on the space $A^{q+1}(X,E)$. Set
\[
a_{lq}(t):=< \langle\phi(t) | \sigma^l(t) \rangle ,\alpha_q>.
\]
Because both $\phi=\phi(t)$ and $\sigma^l(t)$ are analytic in $t$, we have that each $a_{lq}$ is holomorphic in $t$ and so
\[
B(V)=\{t\in B\mid a_{lq}=0, l=1, \cdots, N, q=1, \cdots, d\}
\]
is an analytic subset of $B$.

Next note that by Proposition \ref{prop-holo=harmonic part vanish}
\[
 t\in B(V) \Leftrightarrow V_{t}=V .
\]
So \eqref{B(V)-discription} follows from the fact that $\dim V_{t}=\dim \Image f_t +\dim \ker f_t$. If $V= \mathcal{H}^{0,q}(X,E)$, then $f_t: V_{t} \to H^{0,q}_{\bar{\partial}_{\phi(t)}}(X,E)$ is surjective by Proposition \ref{prop V-E-t}. Hence  \eqref{B'-discription} follows from the extension isomorphism $H^{0,q}_{\bar{\partial}_{\phi}}(X,E) \cong H^{0,q}_{\bar{\partial}_{t}}(X_t,E_t)$, see \eqref{ext-iso for coho}.
\end{proof}

Note that the definition of $B(V)$ depends only on the vector subspace $V\subseteq \mathcal{H}^{0,q}(X,E)$ and the holomorphic tensor bundle $E$. It does not depend on the choices of the basis $\sigma_0^1, \cdots, \sigma_0^N$.

\section{The canonical deformations}
\label{Universal property for the canonical deformation}
Let $\pi: (\mathcal{X}, X)\to (B,0)$ be a deformation of a compact Hermitian manifold $X$ such that for each $t\in B$ the complex structure on $X_t$ is represented by a Beltrami differential $\phi(t)$ and $E$ be a holomorphic tensor bundle on $X$ (see Definition \ref{def-holo-tensor-bundle}). We introduce the following definitions:
\begin{definition}\label{def-deformation-Dol-class}
Given $y\in \ker\bar{\partial}\cap A^{0,q}(X,E)$ and $T\subseteq B$, which is a complex subspace of $B$ containing $0$, a \emph{deformation} of $y$ (w.r.t. $\pi: (\mathcal{X}, X)\to (B,0)$ ) on $T$ is a family of $E$-valued $(0,q)$-forms $\sigma (t)$ such that
\begin{itemize}
  \item[1.] $\sigma (t)$ is holomorphic in $t$;
  \item[2.] $\bar{\partial}_{\phi(t)}\sigma (t) = \bar{\partial}\sigma (t) - \langle\phi (t) | \sigma (t) \rangle =0,~\forall t\in T$;
  \item[3.] $\sigma (0) = y$.
\end{itemize}
If $B$ is smooth and there exists $\sigma (t)\in A^{0,q}(X,E)[[t]]$ (i.e. a formal power series) such that $2.,3.$ are satisfied formally, then we call $\sigma (t)$ a \emph{formal deformation} of $y$.
We say $y$ has (formally) \emph{unobstructed deformation w.r.t. $\pi$} if $B$ is smooth and there exists a (formal) deformation of $y$ on $B$. We say $y$ has (formally) \emph{unobstructed deformation} if for any small deformation $\pi: (\mathcal{X}, X)\to (B,0)$ of $X$ with smooth $B$
there exists a (formal) deformation of $y$ (w.r.t. $\pi$) on $B$. Such $y$ are also called \emph{unobstructed} $E$-valued $(0,q)$-forms.

Given a Dolbeault cohomology class $\alpha\in H_{\bar{\partial}}^{0,q}(X,E)$, a \emph{deformation} of $\alpha$ (w.r.t. $\pi$) on $T$ is defined as a triple $(y,\sigma (t),T)$ which consisting of a representative $y\in\alpha$ and a deformation $\sigma (t)$ of $y$ (w.r.t. $\pi$) on $T$. In particular, $(0,0,B)$ is called the \emph{trivial deformation}. A \emph{formal deformation} of $\alpha$ is a triple $(y,\sigma (t),T)$ where $\sigma (t)$ is a formal deformation of $y\in\alpha$ (w.r.t. $\pi$) on $T$. Two deformations (w.r.t. $\pi$) $(y,\sigma (t),T)$ and $(y',\sigma' (t),T)$ of $\alpha$ on $T$ are \emph{equivalent} if
\[
[\sigma (t) - \sigma' (t)] = 0 \in H^{0,q}_{\bar{\partial}_{\phi(t)}}(X,E),~\forall t\in T.
\]
We say a Dolbeault cohomology class $\alpha\in H_{\bar{\partial}}^{0,q}(X,E)$ has (formally) \emph{unobstructed deformation} if for any small deformation $\pi: (\mathcal{X}, X)\to (B,0)$ of $X$ with smooth $B$ there exists $y\in\alpha$ such that $y$ has (formally) unobstructed deformation w.r.t. $\pi$.
If for any small deformation $\pi: (\mathcal{X}, X)\to (B,0)$ of $X$ with smooth $B$, every Dolbeault cohomology classes in $H_{\bar{\partial}}^{0,q}(X,E)$ have (formally) unobstructed deformation w.r.t. $\pi$, then we say the \emph{deformations of classes in $H_{\bar{\partial}}^{0,q}(X,E)$ are (formally) unobstructed}.
\end{definition}
\begin{remark}We see that the case $q=n$ is trivial by definition since any $E$-valued $(0,n)$-form $\sigma_0$ automatically has a deformation given by $\sigma (t):=\sigma_0$.
\end{remark}

\begin{definition}\label{def-canonical-deformation}
A deformation $\sigma (t)$ of $y\in \ker\bar{\partial}\cap A^{0,q}(X,E)$ on $T$ (w.r.t. $\pi$) is called \emph{canonical} if it satisfies
$\sigma(t) = y + \bar{\partial}^*G \langle\phi(t) | \sigma(t)\rangle$ for any $t\in T$. A deformation $(y,\sigma (t),T)$ of $[y]\in H_{\bar{\partial}}^{0,q}(X,E)$ (w.r.t. $\pi$) is called \emph{canonical} if $\sigma (t)$ is a canonical deformation of $y$.
\end{definition}

\begin{lemma}\label{uniqueness of canonical deformation}
For any small deformation $\pi: (\mathcal{X}, X)\to (B,0)$ and any $y\in \ker\bar{\partial}\cap A^{0,q}(X,E)$, the canonical deformation (if it exists on some $T\subseteq B$ ) of $y$ is unique.
\end{lemma}
\begin{proof}This follows from the same proof as in Proposition \ref{prop-unique-power-seris-solution}.
\end{proof}
Note that if $y$ is harmonic, then the canonical deformation $\sigma (t)$ of $y$ satisfies
\[
\bar{\partial}^*\sigma (t) = 0~\text{and}~\mathcal{H}\sigma (t)=\mathcal{H}y,~ \forall t\in T.
\]
For any $y\in \ker\bar{\partial}\cap A^{0,q}(X,E)$, we have constructed in Theorem \ref{Deformation of Dolbeault cohomology classes: general case} the canonical deformation of $\mathcal{H}y$ on $B(\mathbb{C}\mathcal{H}y)$.

For a given small deformation $\pi: (\mathcal{X}, X)\to (B,0)$ with smooth $B$, we say $y\in \ker\bar{\partial}\cap A^{0,q}(X,E)$ is \emph{canonically unobstructed w.r.t. $\pi$} if its canonical deformation (w.r.t. $\pi$) exists on $B$ and a class $\alpha\in H_{\bar{\partial}}^{0,q}(X,E)$ is \emph{canonically unobstructed w.r.t. $\pi$} if there is a $y\in\alpha$ such that $y$ is canonically unobstructed w.r.t. $\pi$. If every Dolbeault cohomology classes in $H_{\bar{\partial}}^{0,q}(X,E)$ have canonically unobstructed deformation w.r.t. $\pi$, then we say the \emph{deformations of classes in $H_{\bar{\partial}}^{0,q}(X,E)$ are canonically unobstructed w.r.t. $\pi$}. If these statements hold for any small deformation of $X$, we will correspondingly modify the terminology. For example, we say $y\in \ker\bar{\partial}\cap A^{0,q}(X,E)$ is \emph{canonically unobstructed} if for any small deformation $\pi: (\mathcal{X}, X)\to (B,0)$ with smooth $B$ its canonical deformation w.r.t. $\pi$ exists on $B$.

\subsection{The pullback of deformation of Dolbeault cohomology classes}
Given a small deformation $\pi: (\mathcal{X}, X)\to (B,0)$ of $X$ and a holomorphic map $h:(D,s_0)\to(B,0)$, we can always \emph{pullback} $\pi$ to get a new deformation $\varpi: (\mathcal{Y}, Y_{s_0}) \to (D, s_0)$ such that the following diagram is commutative (c.f.~\cite[pp.49]{Bal10}):
\[
\xymatrix{
(\mathcal{Y}, Y_{s_0})  \ar[r]^{\Phi} \ar[d]^{\varpi} & (\mathcal{X}, X) \ar[d]^{\pi} \\
   (D, s_0)  \ar[r]^{h}    & (B,0) , }
\]
where $\mathcal{Y}:=\mathcal{X}\times_B D=\{(u,s)\in \mathcal{X}\times D\mid \pi(u)=h(s)\}$. Given a deformation $\sigma (t)$ of $y\in \ker\bar{\partial}\cap A^{0,q}(X,E)$ on $T\subseteq B$, it is natural to ask whether we can pullback $\sigma (t)$ to get a deformation w.r.t. $\varpi$.

For simplicity, we assume both $B$ and $D$ are smooth (the general case will follow from this). Let $F:X\times B\to \mathcal{X}$ be a smooth trivialization such that the following diagram is commutative:
\[
  \xymatrix{
  X\times B \ar[dr] \ar[r]^-{F}
                & \mathcal{X} \ar[d]^-{\pi}  \\
                & B, }
\]
we define the pullback of $F$ to be a mapping
\begin{align*}
F_{\varpi}:~~~&Y_{s_0}\times D\longrightarrow \mathcal{Y},\\
&(y,s) \mapsto (F(\Phi_{s_0}(y),h(s)),s) ,
\end{align*}
where $\Phi_{s_0}:Y_{s_0}\to X$ is the restriction of $\Phi$. It is straightforward to check that $F_{\varpi}$ is a diffeomorphism such that the following diagram is commutative:
\[
  \xymatrix{
   & D\\
  Y_{s_0}\times D \ar[ur]^-{pr_2} \ar[r]^-{F_{\varpi}} \ar[d]^-{(\Phi_{s_0},h)}  &  \mathcal{Y}  \ar[u]^-{\varpi}  \ar[d]^-{\Phi} \\
 X\times B \ar[r]^-{F} & \mathcal{X},}
\]
where $pr_2:Y_{s_0}\times D\to D$ is the projection mapping with respect to the second factor $D$. We will call $F_{\varpi}$ the \emph{pullback trivialization of $F$}. Note that if $F$ is a transversely holomorphic trivialization for $\pi$ (i.e. $F\mid_{\{x\}\times B}$ is holomorphic for each $x\in X$, c.f. \cite{Cle05,Ma05}) , then $F_{\varpi}$ is also a transversely holomorphic trivialization for $\varpi$.

We have the following
\begin{lemma}\label{lem-pullback}
Let $\pi: (\mathcal{X}, X)\to (B,0)$ be a small deformation with Beltrami differentials $\phi(t)$ (w.r.t. a trivialization $F$) and $\varpi: (\mathcal{Y}, Y_{s_0}) \to (D, s_0)$ the pullback of $\pi$ via a holomorphic map $h:(D,s_0)\to(B,0)$ as in the above commutative diagram. Denote by $\psi(s)$ the Beltrami differentials of $\varpi$ (w.r.t. $F_{\varpi}$, the pullback trivialization of $F$) and $\Phi_{s}:=\Phi\mid_{Y_{s}}$ for each $s\in D$, then we have
\begin{enumerate}
  \item $\Phi_{s_0}^*\phi(h(s))=\psi(s)\in $ for any $s\in D$, where $\Phi_{s_0}^*:A^{0,1}(X,T_{X}^{1,0})\to A^{0,1}(Y_{s_0},T_{Y_{s_0}}^{1,0})$ is induced by $\Phi_{s_0}$;
  \item $\Phi_{s_0}^{*-1}\langle\phi(h(s))| \Phi_{s_0}^{*}=\langle\psi(s)|$, where $\Phi_{s_0}^*:A^{0,\bullet}(X,E)\to A^{0,\bullet}(Y_{s_0},E_{s_0})$ is induced by $\Phi_{s_0}$;
  \item Given a deformation $\sigma (t)$ of $y\in \ker\bar{\partial}\cap A^{0,q}(X,E)$ on $T\subseteq B$,
  \[
  \tau(s):=\Phi_{s_0}^*\sigma (h(s)),\quad \forall s\in h^{-1}(T),
  \]
  is a deformation of $\Phi_{s_0}^*y$ on $h^{-1}{T}\subseteq B$. Furthermore, if $Y_{s_0}$ is equipped with the induced Hermitian metric from $X$ by $\Phi_{s_0}^*$ and $\sigma (t)$ is canonical, the same is true for $\tau(s)$.
\end{enumerate}
\end{lemma}
\begin{proof}It suffice to show $(1)$ since the other two assertions are direct consequences of the first one. Indeed, let $z^1,\cdots, z^n$ be holomorphic coordinates on $X$, then $z'^{\alpha}:=\Phi_{s_0}^*z^\alpha=z^\alpha\circ\Phi_{s_0},~\alpha=1,\cdots, n$, are holomorphic coordinates on $Y_{s_0}$ because $\Phi_{s_0}$ is biholomorphic. Similarly, let $w^1,\cdots, w^n$ be holomorphic coordinates on $X_{h(s)}$, then $w'^{i}:=\Phi_{s}^*w^i=w^i\circ\Phi_{s},~i=1,\cdots, n$, are holomorphic coordinates on $Y_{s}$ because $\Phi_{s}:Y_{s}\to X_{h(s)}$ is biholomorphic. Set $f_t:=F\mid{X\times\{t\}}$ and $f_s':=F_{\varpi}\mid{Y_{s_0}\times\{s\}}$, then the following diagram is commutative:
\[
\xymatrix{
Y_{s_0}  \ar[r]^{f_s'} \ar[d]^{\Phi_{s_0}} & Y_{s} \ar[d]^{\Phi_{s}} \\
   X  \ar[r]^{f_{h(s)} }    & X_{h(s)} , }
\]
which implies that $\Phi_{s_0}^*\circ f_{h(s)}^*w^i=f_s'^*\circ\Phi_{s}^*w^i=f_s'^*w'^i$. The conclusion then follows from the definition of Beltrami differentials (see Section \ref{The extension operator and extension equation}).
\end{proof}
We will call $\tau(s):=\Phi_{s_0}^*\sigma(h(s))$ the \emph{pullback of $\sigma(t)$}. Throughout this paper, pullback deformations are always defined with respect to a pullback trivialization.
\subsection{Properties of the canonical deformations}
\begin{lemma}\label{lem-dbar-exact-canonically-unobstructed}
For any small deformation $\pi: (\mathcal{X}, X)\to (B,0)$ and any $y\in \Image\bar{\partial}\cap A^{0,q}(X,E)$, the canonical deformation of $y$ exists on $B$ and is equivalent to the trivial deformation $(0,0,B)$.
\end{lemma}
\begin{proof}In fact, by Hodge theory
\[
\Image\bar{\partial}=\left(\Image\bar{\partial}\cap(\Image\bar{\partial}_{\phi(t)}+\Image\bar{\partial}^*)\right)
\oplus \Image\bar{\partial}\cap\ker\bar{\partial}_{\phi(t)}^*=\Image\bar{\partial}\cap(\Image\bar{\partial}_{\phi(t)}+\Image\bar{\partial}^*),
\]
where we have used the fact that $\Image\bar{\partial}\cap\ker\bar{\partial}_{\phi(t)}^*=0$ (see the proof of Proposition \ref{partial-bar*-closed-rep of deformed-Dol-coho}). Since $\Image\bar{\partial}_{\phi(t)}\cap\Image\bar{\partial}^*=0$ (see e.g. Remark \ref{rk-kerdbar*-imagedbarphi}) it follows that
\begin{equation}\label{eq-1353}
\Image\bar{\partial}\subseteq \Image\bar{\partial}_{\phi(t)}\oplus\Image\bar{\partial}^*.
\end{equation}
As in the proof of Proposition \ref{prop-kerdbar*-imagedbarphi}, the following map is an isomorphism
\[
\Image\bar{\partial}_{\phi(t)}\longrightarrow \ker\bar{\partial}\cap\left(\Image\bar{\partial}_{\phi(t)}+\Image\bar{\partial}^*\right):
x\longmapsto x-\bar{\partial}^*\bar{\partial}Gx=x-\bar{\partial}^*G\langle\phi(t) | x\rangle,
\]
with its inverse
\[
\ker\bar{\partial}\cap\left(\Image\bar{\partial}_{\phi(t)}+\Image\bar{\partial}^*\right)\longrightarrow \Image\bar{\partial}_{\phi(t)}:
x_0\longmapsto x(t),
\]
where $x(t)$ is the unique solution of the equation $x(t)=x_0+\bar{\partial}^*G\langle\phi(t) | x(t)\rangle$. Now if $y=x_0\in \Image\bar{\partial}\cap A^{0,q}(X,E)$ and $x(t)$ is its canonical deformation, by \eqref{eq-1353} we have $x_0\in \Image\bar{\partial}\subset\ker\bar{\partial}\cap\left(\Image\bar{\partial}_{\phi(t)}+\Image\bar{\partial}^*\right)$, hence $x(t)\in \Image\bar{\partial}_{\phi(t)}$ for any $t\in B$.
\end{proof}

\begin{proposition}\label{prop-nonharmonic-canonical-deformation}
For any small deformation $\pi: (\mathcal{X}, X)\to (B,0)$ and any $\alpha\in H_{\bar{\partial}}^{0,q}(X,E)$, if the canonical deformation of some representative $y\in\alpha$ exists on $T\subseteq B$ then the canonical deformation of any other representative $y'\in\alpha$ also exists on $T\subseteq B$ and their canonical deformations are equivalent.
\end{proposition}
\begin{proof}Let $\sigma (t)$ and $\sigma' (t)$ be the solution of the equations
\[
\sigma(t) = y + \bar{\partial}^*G \langle\phi(t) | \sigma(t)\rangle \quad\text{and}\quad \sigma'(t) = y' + \bar{\partial}^*G \langle\phi(t) | \sigma'(t)\rangle,
\]
then $\sigma(t)-\sigma'(t)=y-y'+ \bar{\partial}^*G \langle\phi(t) | \sigma(t)-\sigma'(t)\rangle$. Since $y-y'\in\Image\bar{\partial}$, it follows from Lemma \ref{lem-dbar-exact-canonically-unobstructed} that $\sigma(t)-\sigma'(t)\in \Image\bar{\partial}_{\phi(t)}$ for any $t\in B$. By our assumption, $\sigma(t)\in \ker\bar{\partial}_{\phi(t)}$ holds for any $t\in T$ which implies $\sigma'(t)\in \ker\bar{\partial}_{\phi(t)}$ for any $t\in T$.
\end{proof}

\begin{proposition}\label{prop semi-uppercontinuous-analyticaly}
Let $\pi: (\mathcal{X}, X)\to (B,0)$ be a deformation of a compact complex manifold $X$ such that for each $t\in B$ the complex structure on $X_t$ is represented by Beltrami differential $\phi(t)$, then the set $\{t\in B\mid \dim H^{q}(X_t,E_t)\geq k\}$ is an analytic subset\footnote{In fact, a more general result of this type holds, see~\cite[pp.\,210]{GR84}.} of $B$ for any nonnegative integer $k$.
\end{proposition}
\begin{proof}It follows from Proposition \ref{prop V-E-t} that
\begin{align*}
&\{t\in B\mid \dim H^{q}(X_t,E_t)\geq k\}\\
=&\{t\in \mathcal{B}\mid \dim V_t/\ker f_t\geq k\}\\
=&\{t\in \mathcal{B}\mid \dim V_t-\dim\left(\ker\bar{\partial}^*\cap\Image\bar{\partial}_{\phi(t)}\right)\geq k\}.
\end{align*}
The conclusion then follows from Proposition \ref{prop-kerdbar*-imagedbarphi}.
\end{proof}
\begin{remark}Note that the Fr\"olicher spectral sequence on $X_t$ degenerates at $E_1$ if and only if $\sum_{p+q=k}\dim H_{\bar{\partial}_t}^{p,q}(X_t)=b_k$ for any $k$ and
\[
\{t\in B\mid\sum_{p+q=k}\dim H_{\bar{\partial}_t}^{p,q}(X_t)= b_k,~\forall~k\}= B\setminus \bigcup_{0\leq k\leq n}\{t\in B\mid\sum_{p+q=k}\dim H_{\bar{\partial}_t}^{p,q}(X_t)\geq b_k+1\}
\]
where $b_k$ is the $k$-th Betti number of $X$. Hence, a direct consequence of Proposition~\ref{prop semi-uppercontinuous-analyticaly} is that the set
\[
T=\{t\in B\mid \text{the Fr\"olicher spectral sequence on}~X_t~\text{degenerates at}~E_1\}
\]
is an analytic open subset (i.e. complement of analytic subset) of $B$.
In particular, if $B$ is a small open disc in the complex plane with $0\in B$ and $T$ is not empty then we have $T=B$ or $T=B\setminus\{0\}$ which corresponds to the fact that degeneration at $E_1$ of the Fr\"olicher spectral sequence is a deformation open property but not a deformation closed property (in the sense of Popovici~\cite{Pop14}). In fact, the Iwasawa manifold provides explicit example for the later phenomenon (see Remark \ref{rk-deformation nonclosed}), see also~\cite{ES93,AK17b}.
\end{remark}

\begin{theorem}\label{thm-2nd-main}
Let $\pi: (\mathcal{X}, X)\to (B,0)$ be a small deformation of a compact Hermitian manifold $X$ with Beltrami differentials $\phi(t)$. For each $q\geq 0$ and $t\in B$, set
\[
v^q_t:=\dim H_{\bar{\partial}}^{0,q}(X,E)-\dim \ker\bar{\partial}_{\phi(t)}\cap\ker\bar{\partial}^*\cap A^{0,q}(X,E) \geq 0,
\]
then we have $(v^{-1}_t:=0)$
\begin{equation}\label{eq-1461}
\dim H_{\bar{\partial}}^{0,q}(X,E)=\dim H_{\bar{\partial}_t}^{0,q}(X_t,E_t)+v^q_t+v^{q-1}_t.
\end{equation}
In particular, $\dim H_{\bar{\partial}_t}^{0,q}(X_t,E_t)$ is independent of $t\in B$ if and only if the canonical deformations of classes in $H_{\bar{\partial}}^{0,q}(X,E)$ and $H_{\bar{\partial}}^{0,q-1}(X,E)$ exist on $B$.
\end{theorem}
\begin{proof}By the second assertion of Proposition \ref{prop-kerdbar*-imagedbarphi} we have
\[
v^{q-1}_t=\dim \Image\bar{\partial}_{\phi(t)}\cap\ker\bar{\partial}^*\cap A^{0,q}(X,E),
\]
which, combined with Proposition \ref{partial-bar*-closed-rep of deformed-Dol-coho} gives $v^q_t+v^{q-1}_t=\dim H_{\bar{\partial}}^{0,q}(X,E)-\dim H_{\bar{\partial}_t}^{0,q}(X_t,E_t)$.

Now it follows from \eqref{eq-1461} that $\dim H_{\bar{\partial}_t}^{0,q}(X_t,E_t)$ is independent of $t\in B$ if and only if $v^q_t=v^{q-1}_t=0$ for any $t\in B$. This is equivalent to the harmonic representatives of classes in $H_{\bar{\partial}}^{0,q}(X,E)$ and $H_{\bar{\partial}}^{0,q-1}(X,E)$ are canonically unobstructed w.r.t. $\pi$ in view of the proof of Proposition \ref{prop V-E-t}. The conclusion then follows from Proposition \ref{prop-nonharmonic-canonical-deformation}.
\end{proof}

We summarize the following properties of the canonical deformation:
\begin{theorem}\label{universal and uniqueness of the canonical deformation}
Let $\pi: (\mathcal{X}, X)\to (B,0)$ be a small deformation of a compact Hermitian manifold $X$ with Beltrami differentials $\phi(t)$.
\begin{itemize}
	\item[$(i)$] Assume $S$ is an analytic subset of $B$ with $0\in S$ and $y\in \ker\bar{\partial}\cap A^{0,q}(X,E)$. Then the canonical deformation of $y$ exists on $B(\mathbb{C}\mathcal{H}y)$ and if the canonical deformation of $y$ exists on $S$ we must have $S\subseteq B(\mathbb{C}\mathcal{H}y)$;
	\item[$(ii)$] For any deformed Dolbeault cohomology class $[u]\in H_{\bar{\partial}_{\phi(t)}}^{0,q}(X,E)$, there exists $\sigma_0\in \mathcal{H}^{0,q}(X,E)$ such that $[u]=[\sigma(t)]$ where $\sigma(t)$ is the canonical deformation of $\sigma_0$;
    \item[$(iii)$] For any class $\alpha\in H_{\bar{\partial}}^{0,q}(X,E)$, if the canonical deformation of some representative $y\in\alpha$ exists on $S\subseteq B$ then the canonical deformation of any other representative $y'\in\alpha$ also exists on $S\subseteq B$ and their canonical deformations are equivalent.
    \item[$(iv)$]  Let $\varpi: (\mathcal{Y}, Y_{s_0}) \to (D, s_0)$ be a pullback of the Kuranishi family $\pi_K: (\mathcal{X},X)\to (\mathcal{B},0)$ with the following commutative diagram:
\[
\xymatrix{
(\mathcal{Y}, Y_{s_0})  \ar[r]^{\Phi} \ar[d]^{\varpi} & (\mathcal{X}, X) \ar[d]^{\pi_K} \\
   (D, s_0)  \ar[r]^{h}    & (\mathcal{B},0) , }
\]
and $\Phi_{s}:=\Phi\mid_{Y_{s}}$ for each $s\in D$. Assume $Y_{s_0}$ is equipped with the induced Hermitian metric so that $\Phi_{s_0}^*:Y_{s_0}\to X$ is an isometry and $\tau(s)$ is the canonical deformation of $\Phi_{s_0}^*y$ on $\tilde{D}\subseteq D$. Then $\tilde{D}\subseteq h^{-1}(\mathcal{B}(\mathbb{C}\mathcal{H}y))$ and $\tau(s)$ is just the pullback of the canonical deformation of $y$ on $\mathcal{B}(\mathbb{C}\mathcal{H}y)$.
\end{itemize}
\end{theorem}
\begin{proof}The first assertion follows from Theorem \ref{Deformation of Dolbeault cohomology classes: general case} and the definition of $B(\mathbb{C}\mathcal{H}y)$. The second assertion follows from Proposition \ref{partial-bar*-closed-rep of deformed-Dol-coho} and Proposition \ref{prop V-E-t}. The third one is just Proposition \ref{prop-nonharmonic-canonical-deformation}. The last one follows from $(i)$, Lemma \ref{uniqueness of canonical deformation} and Lemma \ref{lem-pullback}.
\end{proof}

\begin{corollary}\label{coro-unobstructed classes w.r.t.-varpi}
For $V=\mathcal{H}^{0,q}(X,E)$, we have
\begin{equation}\label{unobstructed classes w.r.t.-varpi-discription}
\{\text{canonically unobstructed classes w.r.t.}~\pi\}=\bigcap_{t\in B}V_{t}
\end{equation}
and
\begin{equation}\label{unobstructed classes w.r.t.-varpi-ineq}
\dim H_{\bar{\partial}_t}^{0,q}(X_t,E_t)+\dim \ker f_t\geq \dim\{\text{canonically unobstructed classes w.r.t.}~\pi\},~~\forall t\in B.
\end{equation}
\end{corollary}
\begin{proof}\eqref{unobstructed classes w.r.t.-varpi-discription} is a direct consequence of the definition of canonical deformations.
\eqref{unobstructed classes w.r.t.-varpi-ineq} then follows from \eqref{unobstructed classes w.r.t.-varpi-discription} and the fact (see Theorem \ref{Deformation of Dolbeault cohomology classes: general case}) that $\dim H_{\bar{\partial}_t}^{0,q}(X_t,E_t)+\dim \ker f_t=\dim V_{t}$, $\forall t\in B$.
\end{proof}

\begin{corollary}\label{coro-deformation invariant implies canonical iso}
If $\dim H_{\bar{\partial}}^{0,q}(X,E)=\dim H_{\bar{\partial}_t}^{0,q}(X_t,E_t)$ for some $t\in B$, then there exists a canonical isomorphism
\[
H_{\bar{\partial}}^{0,q}(X,E)\longrightarrow H_{\bar{\partial}_t}^{0,q}(X_t,E_t):~ [y]~ \mapsto ~ \rho[\sigma(t)]~ ,
\]
where $\sigma (t)$ is the canonical deformation of $y$ and $\rho:H_{\bar{\partial}_{\phi(t)}}^{0,q}(X,E)\to H_{\bar{\partial}_t}^{0,q}(X_t,E_t)$ is the extension isomorphism.
\end{corollary}
\begin{proof}
If $\dim H_{\bar{\partial}}^{0,q}(X,E)=\dim H_{\bar{\partial}_t}^{0,q}(X_t,E_t)$ for some $t\in B$, then by Theorem \ref{thm-2nd-main} we have $v^q_t=v^{q-1}_t=0$. But $v^q_t=0\Leftrightarrow V_t=\mathcal{H}^{0,q}(X,E)$ and $v^{q-1}_t=0\Leftrightarrow\ker f_t\cong\Image\bar{\partial}_{\phi(t)}\cap\ker\bar{\partial}^*\cap A^{0,q}(X,E)=0$. Now the assertion follows from Proposition \ref{prop V-E-t} and the proof of Proposition \ref{prop-nonharmonic-canonical-deformation}.
\end{proof}
\section{Deformations of $(p,q)$-forms}
\label{unobstructed deformations for Kahler mfds and Calabi-Yau mfds}
Recall that the Aeppli cohomology of a complex manifold $X$ is defined as
\[
H^{p,q}_{A}(X) := \frac{\ker \partial\bar{\partial} \cap  A^{p,q}(X)}{ \partial A^{p-1,q}(X)+\bar{\partial}A^{p,q-1}(X) }.
\]
For any $(p,q)\in \mathbb{Z}^2$, consider the following maps
\[
\partial_{A,\bar{\partial}(\ker\partial)}^{p,q}:H_{A}^{p,q}(X)\longrightarrow \frac{\ker\bar{\partial}\cap A^{p+1,q}(X)}{\bar{\partial}\left(\ker\partial\cap A^{p+1,q-1}(X) \right)},\quad \partial_{A,\bar{\partial}}^{p,q}:H_{A}^{p,q}(X)\longrightarrow H_{\bar{\partial}}^{p+1,q}(X),
\]
induced by $\partial$. It is easy to see that $\partial_{A,\bar{\partial}(\ker\partial)}^{p,q}=0$ if and only if for any $\partial \sigma\in \ker\bar{\partial}\cap A^{p+1,q}(X)$ there exists $x\in \ker\partial\cap A^{p+1,q-1}(X)$ such that $\partial \sigma=\bar{\partial}x$. Similarly, $\partial_{A,\bar{\partial}}^{p,q}=0$ if and only if for any $\partial \sigma\in \ker\bar{\partial}\cap A^{p+1,q}(X)$ there exists $x\in A^{p+1,q-1}(X)$ such that $\partial \sigma=\bar{\partial}x$. Note that $\partial_{A,\bar{\partial}}^{p,q}=0$ is also equivalent to the injectivity of $H_{\bar{\partial}}^{p+1,q}(X)\to H_{A}^{p+1,q}(X)$. Similar conditions\footnote{We have been informed by the anonymous referee that by using the structure theory of double complexes one can show $\partial_{A,\bar{\partial}(\ker\partial)}^{p,q}=0$ holds for all $(p,q)\in \mathbb{Z}^2$ is equivalent to the $\partial\bar{\partial}$-lemma on $X$, see \cite{KQ20,PSU21,PSU20,Ste21} for more information.} involving the Bott-Chern cohomology were also considered in \cite[pp.\,688]{Pop19} and \cite{RWZ21,RZ18,AU17,AU16,FY11}, see \cite{Ale19} for a unified discussion.
\begin{theorem}\label{thm-unobstructed-def for p,q-forms}
Let $\pi: (\mathcal{X}, X)\to (B,0)$ be a small deformation of a compact Hermitian manifold $(X,h)$ with Beltrami differentials $\phi(t)$ and such that $B$ is smooth.
\begin{itemize}
  \item[$(i)$] Assume $\partial_{A,\bar{\partial}(\ker\partial)}^{p,q}=0$ and $\partial_{A,\bar{\partial}}^{p-1,q+1}=0$, then the deformations of classes in $H_{\bar{\partial}}^{p,q}(X)$ are canonically unobstructed;
  \item[$(ii)$] Assume $\partial_{A,\bar{\partial}(\ker\partial)}^{p,r}=0$ and $\partial_{A,\bar{\partial}}^{p-1,r+1}=0$ for $r=q,q-1$, then $\dim H_{\bar{\partial}_t}^{p,q}(X_t)$ is independent of $t\in B$.
\end{itemize}
\end{theorem}
\begin{proof}We only need to prove $(i)$ since $(ii)$ follows immediately from $(i)$ and Theorem \ref{thm-2nd-main}.

Let $\sigma_0\in\ker\bar{\partial}\cap A^{p,q}(X)$ and $\sigma(t)=\sum_k\sigma_k$ its canonical deformation w.r.t. $\pi$. It follows from $\partial_{A,\bar{\partial}(\ker\partial)}^{p,q}=0$ that
\[
\partial\sigma_0=0\in \frac{\ker\bar{\partial}\cap A^{p+1,q}(X)}{\bar{\partial}\left(\ker\partial\cap A^{p+1,q-1}(X) \right)}
\]
which means there is a $x_0\in \ker\partial\cap A^{p+1,q-1}(X)$ such that $\partial\sigma_0=\bar{\partial}x_0\Rightarrow i_{\phi_1}\partial\sigma_0=\bar{\partial}i_{\phi_1}x_0\in \Image \bar{\partial}$. Moreover, it follows from $\partial_{A,\bar{\partial}}^{p-1,q+1}=0$ that $\partial i_{\phi_1}\sigma_0\in\ker \bar{\partial}\Rightarrow \partial i_{\phi_1}\sigma_0\in \Image \bar{\partial}$. Hence, $\mathcal{L}_{\phi_1}^{1,0}\sigma_0=\bar{\partial}\sigma_1$ holds. As a result, $\partial\bar{\partial}\sigma_1=\partial i_{\phi_1}\partial\sigma_0 =\partial i_{\phi_1}\bar{\partial}x_0=\partial\bar{\partial}i_{\phi_1}x_0$, it follows from $\partial_{A,\bar{\partial}(\ker\partial)}^{p,q}=0$ that there is a $x_1\in \ker\partial\cap A^{p+1,q-1}(X)$ such that $\partial\sigma_1=\partial i_{\phi_1}x_0+\bar{\partial}x_1$.

Now we prove by induction. Assume the following holds for $1\leq k\leq N$:
\begin{equation}\label{eq-1587}
\left\{
\begin{array}{ll}
\bar{\partial}\sigma_k= \sum_{j=1}^{k}\mathcal{L}_{\phi_j}^{1,0}\sigma_{k-j} , &  \\
\partial\sigma_k=\sum_{j=1}^{k}\partial i_{\phi_j}x_{k-j}+\bar{\partial}x_k, & \\
\end{array} \right.
\end{equation}
where $x_k\in \ker\partial\cap A^{p+1,q-1}(X)$. We need to show \eqref{eq-1587} holds for $k=N+1$. Indeed, it can be proved as in Theorem \ref{deformation of holomorphic sections: special case} (see the proof of \eqref{Existence eq 4}) that $\sum_{j=1}^{N+1}\mathcal{L}_{\phi_j}^{1,0}\sigma_{N+1-j}\in\ker\bar{\partial}$. On the other hand, we have
\begin{align*}
\sum_{j=1}^{N+1}\mathcal{L}_{\phi_j}^{1,0}\sigma_{N+1-j}
&=\sum_{j=1}^{N+1}\left(i_{\phi_j}\partial\sigma_{N+1-j}-\partial i_{\phi_j}\sigma_{N+1-j}\right)\\
&=\sum_{j=1}^{N+1}\left( (\sum_{l=1}^{N+1-j}i_{\phi_j}\partial i_{\phi_l} x_{N+1-j-l})+i_{\phi_j}\bar{\partial}x_{N+1-j}-\partial i_{\phi_j}\sigma_{N+1-j}\right)
\end{align*}
and moreover for each $j$,
\begin{align*}
&i_{\phi_j}\bar{\partial}x_{N+1-j}\\
=&\bar{\partial}i_{\phi_j}x_{N+1-j}-i_{\bar{\partial}\phi_j}x_{N+1-j}\\
=&\bar{\partial}i_{\phi_j}x_{N+1-j}-\frac{1}{2}\sum_{l=1}^{j}(\mathcal{L}_{\phi_l}^{1,0}i_{\phi_{j-l}}-i_{\phi_{j-l}}\mathcal{L}_{\phi_l}^{1,0})x_{N+1-j}\\
=&\bar{\partial}i_{\phi_j}x_{N+1-j}-\frac{1}{2}\sum_{l=1}^{j}(i_{\phi_l}\partial i_{\phi_{j-l}}-\partial i_{\phi_l}i_{\phi_{j-l}}-i_{\phi_{j-l}}i_{\phi_l}\partial +i_{\phi_{j-l}}\partial i_{\phi_l}) x_{N+1-j}\\
=&\bar{\partial}i_{\phi_j}x_{N+1-j}-\frac{1}{2}\sum_{l=1}^{j}(i_{\phi_l}\partial i_{\phi_{j-l}}-\partial i_{\phi_l}i_{\phi_{j-l}} +i_{\phi_{j-l}}\partial i_{\phi_l}) x_{N+1-j}
\end{align*}
which implies
\[
\sum_{j=1}^{N+1}\mathcal{L}_{\phi_j}^{1,0}\sigma_{N+1-j}
=\sum_{j=1}^{N+1}\left(\bar{\partial}i_{\phi_j}x_{N+1-j}+\frac{1}{2}\sum_{l=1}^{j}\partial i_{\phi_l}i_{\phi_{j-l}} \right).
\]
Then it follows from $\partial_{A,\bar{\partial}}^{p-1,q+1}=0$ that there exists a $\sigma_{N+1}$ such that
\[ \sum_{j=1}^{N+1}\mathcal{L}_{\phi_j}^{1,0}\sigma_{N+1-j}=\bar{\partial}\sigma_{N+1}.
\]
Furthermore, note that
\[
\partial\bar{\partial}\sigma_{N+1}=\partial\sum_{j=1}^{N+1}\bar{\partial}i_{\phi_j}x_{N+1-j}.
\]
It follows from $\partial_{A,\bar{\partial}(\ker\partial)}^{p,q}=0$ that $\partial\sigma_{N+1}=\sum_{j=1}^{N+1}\partial i_{\phi_j}x_{N+1-j}+\bar{\partial}x_{N+1}$ for some $x_{N+1}\in \ker\partial\cap A^{p+1,q-1}(X)$.
\end{proof}
\begin{remark}
If $X$ is a compact K\"ahler manifold, then $\sigma_k$ is $\partial$-exact and $\bar{\partial}^*$-exact for each $k>0$. The holomorphic family of $(p,q)$-forms $\sigma(t)$ on $X$, when considered as $(p,q)$-forms on $X\times B$, satisfies
\[
    (\bar{\partial}-\mathcal{L}_{\phi(t)}^{1,0})\sigma(t)|_{X\times B}=0,~~~\text{and}~~~,\partial\sigma(t)|_{X\times B}=0.
\]
This fact was used in an essential way by Clemens~\cite[pp.\,339]{Cle05}.
\end{remark}
\section{Examples of obstructed deformations and the jumping phenomenon}
\label{Examples of obstructed deformations and the jumping phenomena}
Nakamura~\cite{Nak75} classified three-dimensional complex solvable manifolds and computed the Dolbeault cohomology of their small deformations. These provided first examples of the jumping phenomenon. In this section, we analyze these phenomena by using the results obtained in previous sections.

\begin{example}\label{example Case III-(2)}
Case III-(2). Let $G$ be the matrix Lie group defined by
\[
G := \left\{
\left(
\begin{array}{ccc}
 1 & z^1 & z^3 \\
 0 &  1  & z^2 \\
 0 &  0  &  1
\end{array}
\right) \in \mathrm{GL}(3;\mathbb{C}) \mid z^1,\,z^2,\,z^3 \in\mathbb{C} \right\}\cong \mathbb{C}^3,
\]
where the product is the one induced by matrix multiplication. This is usually called the \emph{Heisenberg group}. Consider the discrete subgroup $\Gamma$ defined by
\[
\Gamma := \left\{
\left(
\begin{array}{ccc}
 1 & \omega^1 & \omega^3 \\
 0 &  1  & \omega^2 \\
 0 &  0  &  1
\end{array}
\right) \in G \mid \omega^1,\,\omega^2,\,\omega^3 \in\mathbb{Z}[\sqrt{-1}] \right\},
\]
The quotient $X=G/\Gamma$ is called the \emph{Iwasawa manifold}. A basis of $H^0(X,\Omega^1)$ is given by
\[
\varphi^1 = d z^1,~ \varphi^2 = d z^2,~ \varphi^3 = d z^3-z^1\,d z^2,
\]
and a dual basis $\theta^1, \theta^2, \theta^3\in H^0(X,T_X^{1,0})$ is given by
\[
\theta^1=\frac{\partial}{\partial z^1},~\theta^2=\frac{\partial}{\partial z^2} + z^1\frac{\partial}{\partial z^3},~\theta^3=\frac{\partial}{\partial z^3}.
\]
$X$ is equipped with the Hermitian metric $\sum_{i=1}^3\varphi^i\otimes\bar{\varphi}^i$. The Beltrami differential of the Kuranishi family of $X$ is
\[
\phi(t) = \sum_{i=1}^3\sum_{\lambda=1}^2t_{i\lambda}\theta^i\bar{\varphi}^{\lambda} - D(t)\theta^3\bar{\varphi}^{3},~\text{with}~D(t)=t_{11}t_{22}-t_{21}t_{12},
\]
and the Kuranishi space of $X$ is
\[
\mathcal{B}=\{t=(t_{11}, t_{12}, t_{21}, t_{22}, t_{31}, t_{32})\in \mathbb{C}^6\mid |t_{i\lambda}|<\epsilon, i=1, 2, 3, \lambda=1,2 \},
\]
where $\epsilon>0$ is sufficiently small. Set
\[
\phi_1=\sum_{i=1}^3\sum_{\lambda=1}^2t_{i\lambda}\theta^i\bar{\varphi}^{\lambda},~ \phi_2 = -D(t)\theta^3\bar{\varphi}^{3},
\]
and write the canonical deformation of $\sigma_0\in H_{\bar{\partial}}^{0,q}(X,E)$ by $\sigma(t)=\sum_k \sigma_k(t)$ with each $\sigma_k=\sigma_k(t)$ being the homogeneous term of degree $k$ in $t\in \mathcal{B}$. We need to check that
\begin{equation}\label{obstruction eq 1}
\sum_{j=1}^{k} \langle \phi_j | \sigma_{k-j} \rangle=0\in H_{\bar{\partial}}^{0,q+1}(X,E) , ~~\forall~k>0~ .
\end{equation}

Let us now consider the deformation of classes in
\[
H^{0}(X, \Omega_{X}^{2})=\mathbb{C}\{\varphi^1\wedge\varphi^2, \varphi^2\wedge\varphi^3, \varphi^1\wedge\varphi^3 \}.
\]
First, we compute
\begin{align*}
&\mathcal{L}_{\phi_1}^{1,0}\varphi^1 = \mathcal{L}_{\phi_1}^{1,0}\varphi^2=\mathcal{L}_{\phi_1}^{1,0}\bar{\varphi}^1=\mathcal{L}_{\phi_1}^{1,0}\bar{\varphi}^{2}=\mathcal{L}_{\phi_1}^{1,0}\bar{\varphi}^{3}=0,\\[5pt]
&\mathcal{L}_{\phi_1}^{1,0}\varphi^3 = \sum_{\lambda=1}^2(t_{1\lambda}\varphi^2-t_{2\lambda}\varphi^1)\wedge\bar{\varphi}^{\lambda},~\mathcal{L}_{\phi_2}^{1,0}\varphi^i =\mathcal{L}_{\phi_2}^{1,0}\bar{\varphi}^{i}=0,~i=1,~2,~3.
\end{align*}
and
\begin{align*}
&\mathcal{L}_{\phi_1}^{1,0}(\varphi^1\wedge\varphi^2) = \mathcal{L}_{\phi_1}^{1,0}\varphi^1\wedge\varphi^2 - \varphi^1\wedge\mathcal{L}_{\phi_1}^{1,0}\varphi^2= 0, \\[5pt]
&\mathcal{L}_{\phi_1}^{1,0}(\varphi^2\wedge\varphi^3) = \mathcal{L}_{\phi_1}^{1,0}\varphi^2\wedge\varphi^3 - \varphi^2\wedge\mathcal{L}_{\phi_1}^{1,0}\varphi^3= -t_{2\lambda} \bar{\varphi}^\lambda\wedge\varphi^1\wedge\varphi^2, \\[5pt]
&\mathcal{L}_{\phi_1}^{1,0}(\varphi^1\wedge\varphi^3) = \mathcal{L}_{\phi_1}^{1,0}\varphi^1\wedge\varphi^3 - \varphi^1\wedge\mathcal{L}_{\phi_1}^{1,0}\varphi^3= -t_{1\lambda} \bar{\varphi}^\lambda\wedge\varphi^1\wedge\varphi^2.
\end{align*}
Set $\sigma_0 = a_{12}\varphi^1\wedge\varphi^2 + a_{23}\varphi^2\wedge\varphi^3 + a_{13}\varphi^1\wedge\varphi^3$, then
\[
\mathcal{L}_{\phi_1}^{1,0}\sigma_0= -(t_{2\lambda}a_{23} + t_{1\lambda}a_{13}) \bar{\varphi}^\lambda\wedge\varphi^1\wedge\varphi^2
\]
is exact if and only if $\mathcal{L}_{\phi_1}^{1,0}\sigma_0=0$, i.e.
\begin{equation}\label{a_{23}}
\left\{
\begin{array}{rcl}
t_{21}a_{23} + t_{11}a_{13} &=& 0 \\[5pt]
t_{22}a_{23} + t_{12}a_{13} &=& 0.
\end{array}
\right.
\end{equation}
has solutions for $(a_{23}, a_{13})$, and in this case the (degree $1$ term of) canonical deformation is given by
\[
\sigma_{1}=\bar{\partial}^*G \mathcal{L}_{\phi_1}^{1,0}\sigma_0 =0.
\]
On the other hand,
\[
\mathcal{L}_{\phi_2}^{1,0}\sigma_0 = 0 \Longrightarrow \sigma_{2}=\bar{\partial}^*G (\mathcal{L}_{\phi_2}^{1,0}\sigma_0+\mathcal{L}_{\phi_1}^{1,0}\sigma_1) =0,
\]
and $\phi_k=0,~ k>2$  implies that $\sigma_{k}=0,~ k>2$.

Therefore, for $V=H^{0}(X, \Omega_{X}^{2})$ we have (See Definition \ref{def-V_t})
\[
V_t=\{a_{12}\varphi^{12} + a_{23}\varphi^{23} + a_{13}\varphi^{13} \mid (a_{12}, a_{23}, a_{13})\in \mathbb{C}^3~\text{satisfy}~ \eqref{a_{23}}\}\subseteq H^{0}(X, \Omega_{X}^{2})
\]
where we denote $\varphi^i\wedge\varphi^j$ by $\varphi^{ij}$ for short. Note that $\dim V_t$ is determined by the rank of the matrix
\[
T=
\left(
\begin{array}{ccc}
 0 &  t_{21}  & t_{11} \\
 0 &  t_{22}  &  t_{12}
\end{array}
\right).
\]
Set $h^{p,q}(X_t)=\dim H^{q}(X_t, \Omega_{X_t}^{p})$ and write $(i), (ii), (iii)$ for the three cases when $(t_{11}, t_{12}, t_{21}, t_{22})=0$, $(t_{11}, t_{12}, t_{21}, t_{22})\neq0$ and $D(t)= 0$, $D(t)\neq 0$, respectively. Then by Proposition \ref{prop V-E-t} and Theorem \ref{Deformation of Dolbeault cohomology classes: general case} we have the following
\vspace{12pt}

\renewcommand\arraystretch{1.5}
\begin{table}[!htbp]
\centering
\begin{center}
\begin{tabular}{|c|c|c|c|}
\hline
$t\in \mathcal{B}$ & $\text{rank}~ T$ &  $h^{2,0}(X_t)$ & $H^{0}(X_t, \Omega_{X_t}^{2})$ \\
\hline
$(i)$ & $0$ &  $3$ & $e^{i_{\phi(t)}} \mathbb{C}\{\varphi^{12}, \varphi^{23}, \varphi^{13} \}$  \\
\hline
$(ii)$ & $1$ &  $2$ & $e^{i_{\phi(t)}} \mathbb{C}\{\varphi^{12}, t_{21}\varphi^{13}- t_{11}\varphi^{23}~\text{or}~t_{22}\varphi^{13}- t_{12}\varphi^{23} \}$\\
\hline
$(iii)$ & $2$ &  $1$ & $e^{i_{\phi(t)}} \mathbb{C}\{\varphi^{12} \}$\\
\hline
\end{tabular}
\end{center}
\end{table}
Notice that the canonical deformation of $\varphi^{12}$ on $\mathcal{B}$ is just itself, but it also has other non-canonical deformations such as $(t _{11}+1)\varphi^{12}$.

For the deformation of classes in
\[
H^{1}(X, \Omega_{X}^{1})=\mathbb{C}\{\varphi^{1\bar{1}}, \varphi^{1\bar{2}}, \varphi^{2\bar{1}}, \varphi^{2\bar{2}}, \varphi^{3\bar{1}}, \varphi^{3\bar{2}} \},
\]
we set $\sigma_0 = a_{11}\varphi^{1\bar{1}} +a_{12}\varphi^{1\bar{2}} +a_{21}\varphi^{2\bar{1}} +a_{22}\varphi^{2\bar{2}} +a_{31}\varphi^{3\bar{1}} +a_{32}\varphi^{3\bar{2}}$, then
\begin{align*}
\mathcal{L}_{\phi_1}^{1,0}\sigma_0
&= (a_{31}t_{22} - a_{32}t_{21})\varphi^{1\bar{1}\bar{2}} + (a_{32}t_{11} - a_{31}t_{12})\varphi^{2\bar{1}\bar{2}}\\
&=-\bar{\partial}\left((a_{31}t_{22} - a_{32}t_{21})\varphi^{1\bar{3}}+(a_{32}t_{11} - a_{31}t_{12})\varphi^{2\bar{3}}\right) ,
\end{align*}
and
\[
\sigma_1=\bar{\partial}^*G \mathcal{L}_{\phi_1}^{1,0}\sigma_0 =-(a_{31}t_{22} - a_{32}t_{21})\varphi^{1\bar{3}}-(a_{32}t_{11}- a_{31}t_{12})\varphi^{2\bar{3}},
\]
where in the last equality we have used the fact that $\bar{\partial}^*G \bar{\partial}\varphi^{i\bar{3}}=\varphi^{i\bar{3}}$ since $\varphi^{i\bar{3}}\in \Image \bar{\partial}^*$ for $i=1, 2$. Furthermore,
\[
\mathcal{L}_{\phi_1}^{1,0}\sigma_1=\mathcal{L}_{\phi_2}^{1,0}\sigma_0=\mathcal{L}_{\phi_2}^{1,0}\sigma_1 = 0\Rightarrow \sigma_{k}=0,~ k>1.
\]
We see that classes in $H^{1}(X, \Omega_{X}^{1})$ are all canonically unobstructed. On the other hand, for the deformation of classes in
\[
H^{0}(X, \Omega_{X}^{1})=\mathbb{C}\{\varphi^1, \varphi^2, \varphi^3\}
\]
we set $\sigma_0 = a_1\varphi^1+ a_2\varphi^2+ a_3\varphi^3$, then
\[
\mathcal{L}_{\phi_1}^{1,0}\sigma_0=a_3\sum_{\lambda=1}^2(t_{1\lambda}\varphi^2-t_{2\lambda}\varphi^1)\wedge\bar{\varphi}^{\lambda}
\]
is exact if and only if $a_3t_{11}=a_3t_{12}=a_3t_{21}=a_3t_{22}=0$ and in which case $\sigma_{k}=0$ for all $k>0$.
Therefore, for $E=\Omega_{X}^{1}$, we have the following (in the notation of Theorem \ref{thm-2nd-main})
\vspace{12pt}
\renewcommand\arraystretch{1.5}
\begin{table}[!htbp]
\centering
\begin{center}
\begin{tabular}{|c|c|c|c|c|}
\hline
$t\in \mathcal{B}$ & $h^{1,1}(X)$ & $ v_t^0$ &  $ v_t^1$ & $h^{1,1}(X_t)$ \\
\hline
$(i)$ & $6$ & $0$ & $0$ & $6$  \\
\hline
$(ii)$~\text{and}~$(iii)$ & $6$ & $1$ &  $0$ & $5$\\
\hline
\end{tabular}
\end{center}
\end{table}

Note that $\bar{\partial}_{\phi(t)}\varphi^3$ is a non-canonical deformation of $0$ on $\mathcal{B}$ which is equivalent to the trivial deformation.

Furthermore, we make the following observation: while the canonical deformation of a given Dolbeault cohomology class is unique, it may have inequivalent deformations in general. We can see this from the above example. Consider deformations of classes in $H^{2,3}_{\bar{\partial}}(X)$ where $X$ is the Iwasawa manifold. It is easy to find that $\mathcal{H}^{2,3}_{\bar{\partial}}(X)=\mathbb{C}\{\varphi^{12\overline{123}},\varphi^{23\overline{123}},\varphi^{13\overline{123}}\}$ and we claim
\[
\mathcal{H}^{2,3}_{\bar{\partial}_{\phi(t)}}(X)=
\left\{
\begin{array}{ll}
\mathbb{C}\{\varphi^{12\overline{123}},\varphi^{23\overline{123}},\varphi^{13\overline{123}}\},&\quad (i),\\
\mathbb{C}\{\varphi^{23\overline{123}},\varphi^{13\overline{123}}\},&\quad (ii)~\text{and}~(iii).
\end{array} \right.
\]
In fact, since
\[
(\mathcal{L}_{\phi(t)}^{1,0})^*=\partial^*(i_{\phi(t)})^*-(i_{\phi(t)})^*\partial^*=\partial^*i_{\phi(t)^*}-i_{\phi(t)^*}\partial^*,
\]
where $\phi(t)^*=\sum_{i=1}^3\sum_{\lambda=1}^2\bar{t}_{i\lambda}\bar{\theta}^{\lambda}\varphi^i-\overline{D(t)}\bar{\theta}^3\varphi^{3}$, we may compute
\begin{align*}
\bar{\partial}_{\phi(t)}^*\varphi^{23\overline{123}}&=\bar{\partial}_{\phi(t)}^*\varphi^{13\overline{123}}=0,\\
\bar{\partial}_{\phi(t)}^*\varphi^{12\overline{123}}&=-\sum_{i=1}^2\bar{t}_{i1}\varphi^{i3\overline{23}}+\sum_{i=1}^2\bar{t}_{i2}\varphi^{i3\overline{13}}.
\end{align*}
Our claim is proved. Now we clearly have $\bar{\partial}_{\phi(t)}\varphi^{23\overline{123}}=0$. Hence $\sigma(t):=t_{11}\varphi^{23\overline{123}}$ is a deformation of $[0]\in H^{2,3}_{\bar{\partial}}(X)$ on $\mathcal{B}$. But we notice that this $\sigma(t)$ is not equivalent to the canonical deformation of $0$ (which is identically $0$) because $[\sigma(t)]=[t_{11}\varphi^{23\overline{123}}]\neq 0\in H^{2,3}_{\bar{\partial}_{\phi(t)}}(X)$ if $t\in (ii)$ or $(iii)$.
\end{example}

\begin{remark}\label{rk-deformation nonclosed}
The Kuranishi family of Iwasawa manifold shows that the degenerations at $E_1$ of Fr\"olicher spectral sequence is not a deformation closed property. In fact, let us restrict the Kuranishi family of Iwasawa manifold to the small disc defined by $(t_{12}, t_{21}, t_{31}, t_{32})=0$ and $t_{22}=\varepsilon>0$. Recall that the Fr\"olicher spectral sequence at $E_1$ if and only if the equality holds in the Fr\"olicher inequalities~\cite[pp.\,322]{Dem12}. By the computations of Hodge numbers for the Kuranishi family of the Iwasawa manifold~\cite{Nak75,Ang13}, we see that if $t_{11}=0$ (in class (ii)) the Fr\"olicher spectral sequence on $X_{t}$ does not degenerate at $E_1$ and if $t_{11}\neq 0$ (in class (iii)) the Fr\"olicher spectral sequence on $X_{t}$ degenerates at $E_1$.
\end{remark}
\section{Concluding remarks}
\label{Concluding remarks}
Besides Siu's conjecture about the deformation invariance of plurigenera for K\"ahler manifold, it was asked by Huybrechts~\cite[pp.\,145]{Huy95} whether the dimension of $H^{1}(X, T_{X}^{1,0}\otimes \Omega_{X})$ is a deformation invariant for Calabi-Yau manifolds. It seems that methods are still lacking to solve the extension equation (as a special type of the $\bar{\partial}$-equation) in these cases.

It should be possible to establish a similar theory for more general vector bundles. The most general case would be when the complex structure on $X$ and the vector bundle structure on $E$ varies simultaneously, i.e. deformation of pairs~\cite{Hua95,IM18,CS14,CS18}. Then it seems that we need two Beltrami differentials to capture the deformations, see e.g.~\cite{LT18}. The deformations of Bott-Chern cohomology~\cite{Sch07} will be studied in a subsequent work~\cite{Xia19dBC}. We also believe there is a similar theory for algebraic deformations where the role of Dolbeault cohomology is replaced by Hochschild cohomology~\cite{HG88}.

It is well-known that the Dolbeault cohomology of complex parallelizable manifold may be computed by left invariant forms~\cite{Sak76}. It will be proved in \cite{Xia20g} that given a family of left invariant deformations $\{X_t\}_{t\in B}$ of such manifolds the set of $t$ for which the deformed Dolbeault cohomology may be computed by left invariant forms is an analytic open subset of $B$. This is a refinement of a well-known result of Console-Fino~\cite{CF01}.

It is an interesting question whether the first conclusion of Theorem \ref{thm-unobstructed-def for p,q-forms} still holds if we only assume $\partial_{\bar{\partial},\bar{\partial}}^{p,q}=0=\partial_{\bar{\partial},\bar{\partial}}^{p-1,q+1}$ where
\[
\partial_{\bar{\partial},\bar{\partial}}^{p,q}:H_{\bar{\partial}}^{p,q}(X)\longrightarrow H_{\bar{\partial}}^{p+1,q}(X),
\]
is the natural map induced by $\partial$. This would give a new proof of the deformation invariance of Hodge numbers for compact complex manifolds whose Fr\"olicher spectral sequence degenerates at $E_1$.
\vskip 1\baselineskip \textbf{Acknowledgements.} I would like to thank Prof. Kefeng Liu for his constant encouragement and many useful discussions. Many thanks to Sheng Rao, Quanting Zhao, Guillaume Rond, Xiaokui Yang, Kang Wei, Yang Shen, Shengmao Zhu, Kai Tang, Chunle Huang, Kai Liu and Ruosen Xiong for useful communications. I would also like to thank Prof. Bing-Long Chen for his constant support. I am very grateful to the anonymous referees for their careful reading and for many helpful suggestions.

\bibliographystyle{alpha}

\end{document}